\documentclass[12pt]{article}

\usepackage{amssymb}

\newtheorem{theorem}{$~~~~$ Theorem}[section]

\newtheorem{corollary}[theorem]{$~~~~$ Corollary}
\newtheorem{lemma}[theorem]{$~~~~$ Lemma}

\newtheorem{definition}[theorem]{$~~~~$Definition}

\newtheorem{exx}[theorem]{$~~~~$ Example}
\newtheorem{lemm}[theorem]{$~~~~$ Lemma}

\newtheorem{remm}[theorem]{$~~~~$Remark}

\newtheorem{comm}[theorem]{$~~~~$Comment}

\newtheorem{deff}[theorem]{$~~~~$Definition}

\def\nor1{Normed$\{~2^{ \zzz \theta  \, )} ~$,$~\sqrt{~2^{ \zzz \theta  \, )}}~\}$}

\def\xor2{Normed$\{ ~\sqrt{~2^{ \zzz \theta  \, )}}~,~2~ \} $}

\def\pag2{Page 2}

\def\zzz{~ \sharp ( ~ }

\def\f55{ \normalsize  \baselineskip = 1.8 \normalbaselineskip }

\def\f55{  \baselineskip = 1.1 \normalbaselineskip } 
\def\g55{  \baselineskip = 1.0 \normalbaselineskip } 
\def\s55{ \baselineskip = 1.0 \normalbaselineskip } 

\def\f55{  \baselineskip = 0.7 \normalbaselineskip }

\newcommand{\overx}[1]{\, \overline{ {#1} } \,}

\newcommand{\el}[1]{Line (\ref{#1})}

\newcommand{\eq}[1]{(\ref{#1})}

\begin{document}

\title{About   the 
Chasm Separating the Goals of
Hilbert's Consistency Program From the
Second Incompleteness Theorem}

\def\beq{\begin{equation}}
\def\enq{\end{equation}}

\def\bel{\begin{lemma}}
\def\enl{\end{lemma}}

\def\bec{\begin{corollary}}
\def\enc{\end{corollary}}

\def\bed{\begin{description}}
\def\ennd{\end{description}}
\def\bee{\begin{enumerate}}
\def\ene{\end{enumerate}}

\def\bxbxd{\begin{definition}}
\def\bxbxdd{\begin{definition}}
\def\eedd{\end{definition}}
\def\bxbxdr{\begin{definition} \rm}
\def\bel{\begin{lemma}}
\def\enl{\end{lemma}}
\def\ent{\end{theorem}}

\author{  Dan E.Willard\thanks{This research 
was partially supported
by the NSF Grant CCR  0956495.}}


\date{State University of New York at Albany}

\maketitle

\setcounter{page}{0}
\thispagestyle{empty}

\normalsize

\baselineskip = 1.3\normalbaselineskip

\normalsize

\baselineskip = 1.0 \normalbaselineskip 
\def\bbint{\large \baselineskip = 1.6 \normalbaselineskip } 
\def\bbint{\large \baselineskip = 1.6 \normalbaselineskip }
\def\bbint{\normalsize \baselineskip = 1.3 \normalbaselineskip }


\def\bbint{\normalsize \baselineskip = 1.27 \normalbaselineskip }

\def\bbint{\large \baselineskip = 2.0 \normalbaselineskip }

\def\bbint{\normalsize \baselineskip = 1.25 \normalbaselineskip }
\def\bbina{\normalsize \baselineskip = 1.24 \normalbaselineskip }

\def\bbint{\large \baselineskip = 2.0 \normalbaselineskip }

\def\bbing{ }
\def\bbins{ }
\def\bbinm{ }

\def\bbint{\normalsize \baselineskip = 1.95 \normalbaselineskip }

\def\bbing{ }
\def\bbins{ }
\def\bbinm{ }

\def\bbint{\large \baselineskip = 2.3 \normalbaselineskip } 
\def\bbing{ }
\def\bbins{ }
\def\bbinm{ }

\def\bbint{\normalsize \baselineskip = 1.7 \normalbaselineskip } 

\def\bbint{\large \baselineskip = 2.3 \normalbaselineskip } 
\def\bbinm{ \baselineskip = 1.18 \normalbaselineskip }

\def\bbint{\large \baselineskip = 2.0 \normalbaselineskip } 
\def\bbing{ }
\def\bbins{ }
\def\bbinm{ }
\def\bbinr{ }

\def\bbint{\normalsize \baselineskip = 1.25 \normalbaselineskip }
\def\bbina{\normalsize \baselineskip = 1.24 \normalbaselineskip }
\def\bbinr{ \baselineskip = 1.3 \normalbaselineskip }
\def\bbing{ \baselineskip = 1.28 \normalbaselineskip }
\def\bbins{ \baselineskip = 1.21 \normalbaselineskip }
\def\bbinm{  }

\def\ftl{ \baselineskip = 1.5 \normalbaselineskip }

\bbint

\parskip 5 pt

\noindent

\small

\baselineskip = 1.14 \normalbaselineskip

\parskip 5pt

\baselineskip = 1.2 \normalbaselineskip 



\large
\normalsize

 \baselineskip = 1.2 \normalbaselineskip 


\begin{abstract}
\normalsize
\baselineskip = 1.35 \normalbaselineskip  
We have published several articles about generalizations and 
boundary-case exceptions to the Second Incompleteness Theorem
during the last 25 years. The  current paper will review some of
our prior results and also introduce  an
``enriched'' refinement of  semantic tableaux deduction.
While the Second Incompleteness
Theorem is a strong result, we will 
emphasize its
boundary-case exceptions
are germane to  Global Warming's threat
because our systems can own a simultaneous knowledge about
their own consistency,
together with an understanding of the
$\Pi_1$  implications of
Peano Arithmetic.

\end{abstract}

\bigskip
\bigskip


\bigskip


\bigskip

\small
{\bf Comment:$~$ } The bibliography section of this article contains
citations to all Willard's major
prior
 papers about logic during 1993-2017.

\def\ww22{\normalsize \baselineskip = 1.21\normalbaselineskip \parskip 4 pt}
\def\bb22{\normalsize \baselineskip = 1.19\normalbaselineskip \parskip 4 pt}
\def\zz22z{\normalsize \baselineskip = 1.19 \normalbaselineskip \parskip 3 pt}
\def\xx22{\normalsize \baselineskip = 1.17\normalbaselineskip \parskip 4 pt}
\def\vx22s{\normalsize \baselineskip = 1.16 \normalbaselineskip \parskip 3 pt} 
\def\vv22{\normalsize \baselineskip = 1.17 \normalbaselineskip \parskip 3 pt} 
\def\aa22{\normalsize \baselineskip = 1.15 \normalbaselineskip \parskip 3 pt} 
\def\g55{  \baselineskip = 1.0 \normalbaselineskip } 
\def\s55{ \baselineskip = 1.0 \normalbaselineskip } 
\def\sm55{ \baselineskip = 0.9 \normalbaselineskip }

\vspace*{- 1.0 em}

\def\waw11{\normalsize \baselineskip = 1.72\normalbaselineskip}
\def\waw11{\normalsize \baselineskip = 1.12\normalbaselineskip}
\def\waw11{\normalsize \baselineskip = 1.85\normalbaselineskip}

\def\waw11{\normalsize \baselineskip = 1.45\normalbaselineskip}

\def\waw11{\normalsize \baselineskip = 1.7\normalbaselineskip}

\def\waw11{\normalsize \baselineskip = 1.12\normalbaselineskip}

\def\g55{  \baselineskip = 1.50 \normalbaselineskip } 
\def\s55{ \baselineskip = 1.50 \normalbaselineskip } 
\def\sm55{ \baselineskip = 1.5 \normalbaselineskip }

\def\g55{  \baselineskip = 1.50 \normalbaselineskip } 
\def\s55{ \baselineskip = 1.50 \normalbaselineskip } 
\def\sm55{ \baselineskip = 0.9 \normalbaselineskip }

\def\aa22{\normalsize  \waw11 \parskip 6 pt} 
\def\bb22{\normalsize  \waw11 \parskip 5 pt}
\def\ww22{\normalsize \waw11 \parskip 4 pt}
\def\vv22{\normalsize  \waw11 \parskip 3 pt} 
\def\tt22{\normalsize  \waw11 \parskip 2 pt} 

\def\g55{  \baselineskip = 1.0 \normalbaselineskip } 
\def\b55{  \baselineskip = 1.0 \normalbaselineskip } 
\def\s55{ \baselineskip = 1.0 \normalbaselineskip } 
\def\sm55{ \baselineskip = 0.9 \normalbaselineskip }

\def\mal{ \bf  }
\def\nal{\mathcal}

\def\cvrew{ \baselineskip = 1.6 \normalbaselineskip \parskip 3pt }

\def\ttt2c{ }
\def\tttc{ }

\def\tttc{\tiny \baselineskip = 0.8 \normalbaselineskip  \parskip 0pt }
\def\ttt2c{\tiny \baselineskip = 0.7 \normalbaselineskip  \parskip 0pt }
\def\tttc{ \baselineskip = 2.1 \normalbaselineskip  \parskip 5pt }
\def\ttt2c{ \baselineskip = 2.1 \normalbaselineskip  \parskip 5pt }

\def\tttc{ \baselineskip = 1.15 \normalbaselineskip  \parskip 5pt }
\def\ttt2c{ \baselineskip = 1.15 \normalbaselineskip  \parskip 5pt }

\def\tttc{ \baselineskip = 1.12 \normalbaselineskip  \parskip 4pt }
\def\ttt2c{ \baselineskip = 1.12 \normalbaselineskip  \parskip 4pt }

\def\tttc{ \baselineskip = 1.14 \normalbaselineskip  \parskip 3pt }
\def\ttt2c{ \baselineskip = 1.14 \normalbaselineskip  \parskip 4pt }

\def\cvt{ \baselineskip = 0.98 \normalbaselineskip }
\def\cv9{ \baselineskip = 0.99 \normalbaselineskip }
\def\cvs{ \baselineskip = 1.0 \normalbaselineskip }
\def\cvl{ \baselineskip = 1.0 \normalbaselineskip }
\def\cvh{ \baselineskip = 1.03 \normalbaselineskip }
\def\cvg{ \baselineskip = 1.00 \normalbaselineskip }

\def\cvt{ \baselineskip = 1.6 \normalbaselineskip }
\def\cv9{ \baselineskip = 1.6 \normalbaselineskip }
\def\cvs{ \baselineskip = 1.6 \normalbaselineskip }
\def\cvl{ \baselineskip = 1.6 \normalbaselineskip }
\def\cvh{ \baselineskip = 1.6 \normalbaselineskip }
\def\cvg{ \baselineskip = 1.6 \normalbaselineskip }
\def\cvb{ \baselineskip = 1.6 \normalbaselineskip }
\def\cvnew{ \baselineskip = 1.6 \normalbaselineskip }
\def\cvmew{ \baselineskip = 1.6 \normalbaselineskip }
\def\cvwew{ \baselineskip = 1.6 \normalbaselineskip \parskip 5pt }
\def\cvrew{ \baselineskip = 1.6 \normalbaselineskip \parskip 3pt }

\def\cvt{ \baselineskip = 1.22 \normalbaselineskip }
\def\cv9{ \baselineskip = 1.22 \normalbaselineskip }
\def\cvs{ \baselineskip = 1.22 \normalbaselineskip }
\def\cvl{ \baselineskip = 1.22 \normalbaselineskip }
\def\cvh{ \baselineskip = 1.22 \normalbaselineskip }
\def\cvg{ \baselineskip = 1.22 \normalbaselineskip }
\def\cvb{ \baselineskip = 1.22 \normalbaselineskip }
\def\cvnew{ \baselineskip = 1.4 \normalbaselineskip }
\def\cvmew{ \baselineskip = 1.35 \normalbaselineskip }
\def\cvwew{ \baselineskip = 1.4 \normalbaselineskip \parskip 5pt }
\def\cvrew{ \baselineskip = 1.22 \normalbaselineskip \parskip 3pt }

\def\cvt{ }
\def\cv9{ }
\def\cvs{ }
\def\cvl{ }
\def\cvh{ }
\def\cvg{ }
\def\cvb{ }
\def\cvnew{ } 
\def\cvmew{ }
\def\cvwew{ }
\def\cvrew{ }

\def\fend{ 

\medskip -------------------------------------------------------}

\def\g55{  \baselineskip = 1.0 \normalbaselineskip } 
\def\s55{ \baselineskip = 1.0 \normalbaselineskip } 
\def\sm55{ \baselineskip = 1.0 \normalbaselineskip } 
\def\h55{  \baselineskip = 1.08 \normalbaselineskip } 
\def\b55{  \baselineskip = 1.1 \normalbaselineskip } 

\normalsize

\baselineskip = 1.85 \normalbaselineskip





\parskip 2pt

\vspace*{- 1.0 em}



\baselineskip = 1.04 \normalbaselineskip 
\parskip 2pt

\baselineskip = 0.96 \normalbaselineskip 

%
\baselineskip = 2.16 \normalbaselineskip 
\baselineskip = 2.3 \normalbaselineskip 

\baselineskip = 0.95 \normalbaselineskip 
\baselineskip = 0.88 \normalbaselineskip 
\parskip 0pt

\noindent

%
%
%
%

\newpage

\def\gvs{ \normalsize \baselineskip = 1.4 \normalbaselineskip  \parskip    5pt}
\def\gvs{ \normalsize \baselineskip = 1.44 \normalbaselineskip  \parskip    5pt}
\def\gvs{ \large \baselineskip = 1.44 \normalbaselineskip  \parskip    5pt}
\def\gvs{ \normalsize \baselineskip = 1.44 \normalbaselineskip  \parskip    5pt}\def\gvs{ \normalsize \baselineskip = 1.74 \normalbaselineskip  \parskip    5pt}
\def\gvs{ \normalsize \baselineskip = 1.44 \normalbaselineskip  \parskip 5pt}

\def\gvs{   \baselineskip = 1.74 \normalbaselineskip  \parskip    5pt}

\def\gvs{ \normalsize \baselineskip = 1.44 \normalbaselineskip  \parskip 5pt}
\def\gvs{ \large \baselineskip = 2.0 \normalbaselineskip  \parskip 5pt}
\def\gvs{ \Large \baselineskip = 2.0 \normalbaselineskip  \parskip 5pt}
\def\gvs{ \normalsize \baselineskip = 2.44 \normalbaselineskip  \parskip 5pt}
\def\gvs{ \normalsize \baselineskip = 2.04 \normalbaselineskip  \parskip 5pt}
\def\gvs{ \normalsize \baselineskip = 2.64 \normalbaselineskip  \parskip 5pt}
\def\gvs{ \Large \baselineskip = 1.6 \normalbaselineskip  \parskip 5pt}

\gvs

\footnotesize

\def\gvs{ }

\normalsize \baselineskip = 0.98 \normalbaselineskip
\normalsize \baselineskip = 1.0 \normalbaselineskip
\normalsize \baselineskip = 1.01 \normalbaselineskip

\def\gvs{ \normalsize \baselineskip = 1.25 \normalbaselineskip  \parskip 4pt}

\def\gvs{ \Large \baselineskip = 1.6  \normalbaselineskip  \parskip 6pt}
\def\gvs{ \normalsize \baselineskip = 1.6  \normalbaselineskip  \parskip 6pt}
\def\gvs{ \large \baselineskip = 1.6  \normalbaselineskip  \parskip 6pt}

\def\gvs{ \normalsize \baselineskip = 1.227 \normalbaselineskip  \parskip 3pt}
\def\gvs{ \large \baselineskip = 1.8  \normalbaselineskip  \parskip 6pt}

\def\gvs{ \normalsize \baselineskip = 1.5 \normalbaselineskip  \parskip 3pt}

\def\gvs{ \large \baselineskip = 2.1  \normalbaselineskip  \parskip 6pt}

\def\gvs{ \normalsize \baselineskip = 2.1  \normalbaselineskip  \parskip 6pt}


 \def\gvs{ \normalsize \baselineskip = 1.227 \normalbaselineskip  \parskip 3pt}

 \def\gvs{ \large  \baselineskip = 1.6 \normalbaselineskip  \parskip 5pt}

\def\gvs{ \Large  \baselineskip = 1.8 \normalbaselineskip  \parskip 5pt}
\def\gvs{ \LARGE  \baselineskip = 1.8 \normalbaselineskip  \parskip 5pt}
\def\gvs{ \normalsize  \baselineskip = 2.0 \normalbaselineskip  \parskip 5pt}

\def\gvs{ \Large  \baselineskip = 2.0 \normalbaselineskip  \parskip 5pt}

\def\gvs{ \large  \baselineskip = 2.2 \normalbaselineskip  \parskip 5pt}

\def\gvs{ \normalsize \baselineskip = 2.4  \normalbaselineskip  \parskip 6pt}

\def\gvs{ \normalsize \baselineskip = 2.6  \normalbaselineskip  \parskip 6pt}
\def\gvs{ \normalsize \baselineskip = 2.2  \normalbaselineskip  \parskip 6pt}
\def\gvs{ \normalsize \baselineskip = 1.8  \normalbaselineskip  \parskip 5pt}
\def\tttc{ }
\def\ttt2c{ }

\def\gv2{ \normalsize \baselineskip = 1.30  \normalbaselineskip  \parskip 3pt}

\def\gvs{ }

\def\gvs{ \normalsize \baselineskip = 2.1 \normalbaselineskip  \parskip 7pt}
\def\gvs{ \normalsize \baselineskip = 1.8 \normalbaselineskip  \parskip    7pt}

 \def\gvs{ \large \baselineskip = 1.7  \normalbaselineskip  \parskip 9pt}
\def\gvs{ \normalsize \baselineskip = 2.0  \normalbaselineskip  \parskip 9pt}

\def\gv2{ \normalsize \baselineskip = 1.30  \normalbaselineskip  \parskip 3pt}

\def\gvs{ \large \baselineskip = 1.7  \normalbaselineskip  \parskip 5pt}

\def\gvs{ \normalsize \baselineskip = 2.0  \normalbaselineskip  \parskip 8pt}
\def\gvs{ \large \baselineskip = 2.0  \normalbaselineskip  \parskip 8pt}


\def\gvx{ \large \baselineskip = 1.6 \normalbaselineskip  \parskip    3pt}

\def\gvs{ \Large \baselineskip = 1.8 \normalbaselineskip  \parskip   5pt}

\def\gvs{ \LARGE \baselineskip = 2.0 \normalbaselineskip  \parskip   5pt}

\def\gvs{ \normalsize \baselineskip = 2.0 \normalbaselineskip  \parskip   5pt}
\def\gvs{ \normalsize \baselineskip = 1.4 \normalbaselineskip  \parskip   5pt}

 \def\gvs{ \large \baselineskip = 2.0 \normalbaselineskip  \parskip   5pt}
 \def\gvs{ \Large \baselineskip = 2.0 \normalbaselineskip  \parskip   5pt}
\def\gvs{ \normalsize \baselineskip = 2.0 \normalbaselineskip  \parskip   5pt}
\def\gvs{ \normalsize \baselineskip = 2.5 \normalbaselineskip  \parskip   5pt}

\def\gvs{ \normalsize \baselineskip = 1.4 \normalbaselineskip  \parskip   5pt}

\section{Introduction }

\label{nnn1}
\gvs
\tttc

The existence of a deep chasm between 
the  goals of 
Hilbert's consistency program and the implications of
G\"{o}del's
Second Incompleteness Theorem was
  immediately
apparent
when G\"{o}del
announced his famous
millennial discovery \cite{Go31}.

Interestingly,
neither  G\"{o}del (in 1931) 
nor Hilbert (during the remainder
of his life)
dismissed the
 existence of  possible compromise
solutions, whereby a
{\it  fragment} of the goals of the
Consistency Program could remain intact.
For instance, Hilbert never withdrew
his statement $*$ 
for justifying
his 
 Consistency 
 Program in
 \cite{Hil26}:
\begin{quote}
\small
\baselineskip = 1.2 \normalbaselineskip
\ttt2c 
$*~$
{\it ``
Let us admit that the situation in which we presently
find ourselves with respect to paradoxes is in the long
run intolerable. Just think: in mathematics, this paragon of
reliability and truth, the very notions and inferences,
as everyone learns, teaches, and uses them, lead to absurdities.
And 
where 
else 
would 
reliability and truth be found 
if  even mathematical thinking fails?''}
\end{quote}
%
Indeed, 
the motto of Hilbert's  Consistency Program
({\it ``Wir m\"{u}ssen wissen---Wir werden wissen''} )  
was engraved
onto Hilbert's tombstone \footnote{ English translation: {\it
``We must know: We will know.''}}.

Also, G\"{o}del was
quite  cautious (especially during the early 1930's) not
to speculate
 whether all facets of Hilbert's Consistency program
would come to a termination.
He thus inserted the following
cautious 
caveat into
his
famous 1931 
millennial 
paper  \cite{Go31}:
\begin{quote}
\small
\baselineskip = 1.2 \normalbaselineskip
\ttt2c 
$~**~~$ 
{\it ``It must be expressly noted that
Theorem XI''}
(e.g. the Second Incompleteness Theorem) 
{\it ``represents no contradiction of the formalistic
standpoint of Hilbert. For this standpoint
presupposes only the existence of a consistency
proof by finite means, and {\it there might
conceivably be finite proofs} which cannot
be stated in P or in ... ''}
\end{quote}
Several 
biographies
of 
 G\"{o}del
\cite{Da97,Go5,Yo5}
have noted that
 G\"{o}del's 
 intention (prior to 1930)
was
to
establish
Hilbert's proposed objectives, before
he formalized
his famous result 
that led
in the opposite direction.
For instance,
 Yourgrau's
biography
 \cite{Yo5}
 of G\"{o}del  does record
 how
von Neumann, 
found it necessary during the early 1930's to
{\it ``argue 
against G\"{o}del 
himself''}
 about the definitive 
 termination of Hilbert's
consistency program,  
which
{\it ``for several years''} after \cite{Go31}'s publication,
G\"{o}del 
{\it ``was cautious not to prejudge''}.
It is known that
 G\"{o}del 
had
hinted that
the Second Incompleteness
Theorem was more significant 
during a 1933 Vienna
lecture  \cite{Go33}. (G\"{o}del has been recorded
as having
embraced the Second Incompleteness Theorem
quite
 broadly 
during  the mid-1930's
after he learned
about Turing's work
\cite{Tur36}.)
Yet despite  this endorsement,
 Gerald Sacks  \cite{YouSa14}
has recalled G\"{o}del
telling him that some type of 
revival of Hilbert's Consistency
Program
would be feasible (see footnote 
\footnote{\label{f2}
Some quotes from Sacks's
YouTube talk
\cite{YouSa14} are that G\"{o}del
 {\it ``
 did not  think''}
the objectives of Hilbert's Consistency Program 
{\it ``were erased''} 
by
the Incompleteness Theorem, and
G\"{o}del believed (according to Sacks) 
 it left
  Hilbert's program 
{\it `` very much alive and
even more interesting than it initially was''}. 
} 
for the exact details
that Sacks has provided.)

The research that has followed
 G\"{o}del's  seminal  1931 discovery has 
mainly focused
on studying generalizations of the Second Incompleteness
Theorem
(instead of also examining its
 boundary-case
exceptions). Many of  these generalizations 
of the Second Incompleteness Theorem 
\cite{AZ1,BS76,BI95,Fe60,Fr79a,Ha11,HP91,KT74,Pa71,PD83,Pa72,Pu85,Pu96,So94,Sv7,Vi5,WP87,ww1,ww2,wwlogos,wwapal}  
are quite beautiful. The author of this paper is
            especially impressed by a generalization of the
Second Incompleteness Effect, arrived at by the 
combined work of Pudl\'{a}k and  Solovay
(abetted by
the research of
 Nelson and Wilkie-Paris)
\cite{Ne86,Pu85,So94,WP87}.
These results, which 
  also have been more recently
discussed 
 in \cite{BI95,Ha7,Sv7,ww1},
have noted the Second Incompleteness Theorem
does not require the presence of the Principle of Induction
to apply to most formalisms that use a Hilbert-style form of 
deduction. (The 
next chapter will offer a detailed summary of this 
important generalization
 of
the Second Incompleteness Theorem in
its
 Remark \ref{nremm-2.5}.)

Our research, during the last 25 years has 
had a different focus,
exploring 
Boundary-Case exceptions to the Second Incompleteness Theorem
more intensively than its generalizations.
It would be natural for many readers to ask why such
exceptions should also be studied, so intensively?

The reason is that while generalizations of the Second Incompleteness
Theorem are very pure form a mathematical standpoint, it must
not be forgotten that   Mankind's survival  in the future  will
require developing formalisms that own enough ingenuity to
solve a variety of pressing ecological problems, such as Global 
Warming, in a satisfactory manner.

More specifically, the solution of
mundane problems, threatening human survival, do not require 
use of a formalism 
producing short proofs of
the existence of large integers, whose binary encodings
employ more digits than the
number of
       atoms in the universe.
It is, however, vital
for mature logical formalisms
 to appreciate their own consistency,
{\it in at least a fragmentary sense,} when they
reason about the implications of their own reasoning.
(Otherwise, a Thinking Being, whether computerized or
human, would not be able to explain to itself fully why it is 
of foundational importance to study its own thinking process,
as a fundamental problem-solving facility.)  

In particular, there is no doubt that
a branch of mathematics, that makes it difficult to 
manufacture {\it abbreviated} proofs of the 
existence of numbers as large  as a
google-plex (e.g. $2^{2^{100}}~$ ),
does
 fall 
short of
the Utopian ideals that the intellectual community 
wishes for Mathematics.
We will argue, however, that the striking
engineering needs that
confront
 modern Mankind requires the evolution
of alternate formalisms, however {\it theoretically weak}, for an
adequate result to be obtained for many
 pressing issues.

In other words,
we will contend
that
Hilbert and G\"{o}del
were 
essentially
correct when their statements $*$ and $**$ suggested
that a nihilistic approach, which ignores the engineering-style
capacities of weaker formalisms, that own a fragmented conception of
their own consistency, 
has serious short-comings.
This is because the
dangers of
Global Warming and other 
imminent
threats that endanger 21-st century
Mankind are too serious for logicians to entertain using
anything less than a formalism, which possesses at least a 
fragmentary conception of its own consistency.

Also, we will discuss an addition to our prior research about
self-justifying logics, called {\bf Deductive Enrichment},
which should convince skeptical readers that our formalisms do
indeed have practical value. 
Especially 
within a
special
 context where modern computers can perform
arithmetic operations with more than a
billion times the speed
of a human brain, we will argue our boundary-case exceptions to
the prior century's Second Incompleteness Theorem have noteworthy
pragmatic significance.

As the reader examines this paper, it should be kept in
mind that all our self-justifying axiom systems (since 1993)
contain an ability to prove analogs of all the
$\Pi_1$ theorems of Peano Arithmetic under a 
slightly revised language (such as 
\textsection \ref{s3}'s 
$~L^*$ formalism).
This
fact is non-trivial because an axiom system
that recognizes its own consistency will contain 
little pragmatic 
significance,
 if it does not maintain an
ability to prove all the
quite central
$\Pi_1$ theorems of Peano Arithmetic.

What will make our formalisms tempting
in the current article is 
the  new notion of ``Deductive Enrichment''. It will allow
a formalism to
maintain a simultaneous knowledge about its own consistency
together with 
a recognition about
the
 truth behind the 
$\Pi_1$ theorems of Peano Arithmetic.

In particular, we do not dispute that our formalism will fall
short of the Utopian ideals for mathematics when it is
unable to produce a  {\it brief} proof for the existence of
large numbers, such as a 
google-plex (e.g. $2^{2^{100}}~$ ).
We do, however, claim our formalism contains 
some
pragmatic value
when it
owns a simultaneous knowledge about its own consistency
and the truth behind the 
$\Pi_1$ theorems of Peano Arithmetic.
(Moreover, we
do
 encourage the reader to 
examine, again, G\"{o}del's caution about the
Second Incompleteness Theorem,
which was observed by  Sacks in \cite{YouSa14} and has also been
 summarized by us in
Footnote \ref{f2}.)

\section{General Notation and Literature Survey}

\label{nnn2}

Let us 
call an
ordered pair $(\alpha,d)$ a
    {\bf ``Generalized Arithmetic''}
when  its 
first and second 
components 
are 
defined 
as 
follows:
\bee
\item
The {\bf Axiom Basis} ``$~\alpha~$'' 
of a  Generalized Arithmetic
will be defined as
the set of 
 proper axioms 
it employs.
\item
The second component ``$~d~$'' of a Generalized Arithmetic
will be  defined as
the 
{\it combination} of its formal rules of inference
and 
the
 logical axioms ``$~L_d~$'' it employs.
This second component, ``$~d~$'' of a Generalized Arithmetic
will be called its  {\bf Deductive Apparatus}. 
\ene


\begin{exx}
\label{nex-2.1}
\rm
This notation 
allows us to
 conveniently separate  the logical axioms
$~L_d~,~$ associated with  $(  \alpha  , d  )~$, from 
its
``basis axioms'' $\, \alpha \,$.
It also allows one to  compare
the various
deductive apparatus techniques
that  
have
appeared in the literature.
For instance,
the
  $~d_E~$ apparatus,
introduced
 in 
\textsection
 2.4 of  Enderton's textbook \cite{End}, 
has
 used only  modus ponens
as a rule of inference,
combined with a 
complicated
4-part  schema of logical axioms.
This differs from
the  $~d_M~$ ,  $~d_H~$ and  $~d_F~$ approaches of
Mendelson \cite{Mend}, 
H\'{a}jek-Pudl\'{a}k \cite{HP91}
and Fitting \cite{Fit}.
The former two textbooks
employ a simpler set of logical axioms
than $\, d_E \,$,
but they require
two rules of inference
(modus ponens and generalization).
The  $~d_F~$ apparatus, appearing in Fitting's textbook \cite{Fit},
as well as its predecessor due to Smullyan \cite{Smul},
actually  employ {\it no logical axioms.}
Instead, Fitting and Smullyan rely upon a
``tableaux style'' method for generating a 
consequently
larger number of
rules of
inference.
\end{exx}

\begin{deff}
\label{ndef-2.2}
\rm
Let
$~\alpha~$ once again 
denote an axiom basis, 
and $~d~$ 
designate
 a
deduction apparatus.
Then 
the  ordered pair
 $~(  \alpha  , d  )$
will
be called  a {\bf Self Justifying}
configuration
 when:
\begin{description}
  \item[  i   ] one of  $~(  \alpha  , d  )$'s  theorems
(or possibly one of $\alpha$'s axioms)
do
state that the deduction method $ \, d, \, $ applied to the
basis
system $ \, \alpha, \, $ 
produces a consistent set of theorems, and
\item[  ii   ]
     the axiom system $ \, \alpha  \, $ is in fact consistent.
\end{description}
\end{deff}

\begin{exx}
\label{nex-2.3}
\rm
Using 
Definition \ref{ndef-2.2}'s 
 notation, our
prior
 research  
in
\cite{ww93,ww1,ww5,wwapal,ww9,ww14}
has
developed
arithmetics
$~(  \alpha  , d  )$
that 
 were
``Self Justifying''.
It 
also 
proved
the Second Incompleteness Theorem 
implies specific
limits beyond which 
self-justifying
formalisms cannot transgress.
For any  $\,(\alpha,d) \,$, 
it is 
almost trivial 
to construct a 
system $ \, \alpha^d \, \supseteq  \,  \alpha  \, $
 that  satisfies
the
Part-i 
condition
(in an isolated context {\it where the Part-ii condition is
 not also
satisfied}).
For instance,  $ \, \alpha^d \, $  could
consist of all of $~\alpha \,$'s axioms plus 
an added {\bf $\,$``SelfRef$(\alpha,d)$''$\,$} sentence,
defined as stating:
\begin{quote} 
$\oplus~~~$ 
There is no proof 
(using 
$d$'s deduction method)
of  $0=1$
from the  {\it union}
 of
the
 axiom system $\, \alpha \, $
with {\it this}
sentence  ``SelfRef$(\alpha,d) \,$'' (looking at itself).
\end{quote}
Kleene 
discussed
 in
\cite{Kl38} 
how
to
encode rough
 analogs of 
the above
 {\bf $\,$``I Am Consistent''} 
axiom
statement.
Each of
Kleene, 
Rogers and Jeroslow 
 \cite{Kl38,Ro67,Je71},
however,
 emphasized
that
$\alpha ^d$ 
may
be inconsistent
(e.g.  violating Part-ii of   self-justification's
definition),
{\it despite SelfRef$(\alpha,d)$'s 
formalized
assertion.}
This is because if the 
 pair $(\alpha,d)$ is too strong
then a
quite conventional
G\"{o}del-style diagonalization argument can
be applied to the axiom basis of
$\alpha^d~~=~~ \alpha \, + \, $ SelfRef$(\alpha,d), ~$
where the added presence of the statement 
SelfRef$(\alpha,d)$ 
will cause this extended version of 
$\, \alpha\,$, ironically,
 to
 become automatically inconsistent.
Thus, the
encoding of
``SelfRef$(\alpha,d)$'' is relatively easy,
via an application of the Fixed Point Theorem,
but this sentence
is often,
ironically, 
entirely
useless! 
\end{exx}


\begin{deff}
\label{ndef-2.4}
\rm
Let
 $Add(x,y,z)$ and    $Mult(x,y,z)$ 
denote two 3-way predicate symbols specifying
that $x+y=z$ and $x*y=z$ (under 
 $\Pi_1$ styled-encodings for the
associative, commutative, identity and distributive principles
 using these two 3-way predicate symbols).
Let  $\, \alpha \,$ denote what the first paragraph
of this section 
had
called an ``axiom basis''.
We will 
then
say that $~\alpha~$
{\bf recognizes} successor,  addition  and multiplication
as {\bf Total Functions} iff 
it can  prove
\eq{totdefxs} - \eq{totdefxm}
as theorems.
 {
{
\beq 
\label{totdefxs}
\forall x ~ \exists z ~~~Add(x,1,z)~~
\enq
\beq 
\label{totdefxa}
\forall x ~\forall y~ \exists z ~~~Add(x,y,z)~~
\enq
\beq 
\label{totdefxm}
\forall x ~\forall y ~\exists z ~~~Mult(x,y,z)~
\enq }
}
Also, an axiom basis
$\alpha$ will be  called 
{\bf Type-M} iff it includes
\eq{totdefxs} - \eq{totdefxm}
as theorems, {\bf Type-A} if it includes 
only \eq{totdefxs} and \eq{totdefxa} as theorems,
and {\bf Type-S} if it contains
only \eq{totdefxs} as a
 theorem. 
Also,
$\alpha$ 
is
called 
{\bf Type-NS} iff it  can prove
none of these theorems.
\end{deff}

\bigskip

\begin{remm}
\label{nremm-2.5}
\rm
The separation of basis axiom systems into the four 
categories of Type-NS, Type-S, Type-A and Type-M systems
enables us to nicely summarize the prior literature
about generalizations and boundary-case exceptions
for the Second Incompleteness Theorem. This is because:
\bed
\item[   $~~~~$a.$~~$  ]
The combined research of Pudl\'{a}k, Solovay, Nelson and Wilkie-Paris
\cite{Ne86,Pu85,So94,WP87},
as is formalized by Theorem $\, ++ \,$,
implies no
natural Type-S  generalized arithmetic $(\alpha,d)$
can recognize  its own  consistency
when $d$ is one of 
Example \ref{nex-2.1}'s three
examples of
  Hilbert-style 
deduction operators of 
$\, d_E \,$, $\, d_H \,$  
or
$\, d_M ~~$. In particular,
it establishes the following result:
\begin{quote}
\normalsize \baselineskip = 1.2 \normalbaselineskip 
{\bf ++ }
{\it 
(Solovay's  
modification
\cite{So94}
of Pudl\'{a}k \cite{Pu85}'s formalism 
using some of 
Nelson and Wilkie-Paris \cite{Ne86,WP87}'s
methods)} :
Let 
$ \, (\alpha,d) \, $ 
denote 
a generalized arithmetic
supporting
the
\el{totdefxs}'s
Type-S statement
and 
assuring
the successor operation
will
satisfy
both 
 $  \,   x'     \neq 0     $ and
$     x'     =     y' \Leftrightarrow x=y $.
$~$Then
$ \, (\alpha,d) \, $  
cannot verify its own
consistency
whenever
simultaneously
 $d$ is
a Hilbert-style 
deductive
apparatus and
$~\alpha~$
 treats addition and multiplication
as 3-way relations, 
satisfying 
their usual 
associative, commutative 
 distributive 
and identity 
axioms.
\end{quote}
Essentially, Solovay \cite{So94} 
privately communicated 
to us 
in 1994
an analog of theorem $++$.
Many authors
have noted Solovay
 has 
been
reluctant to publish
his 
nice 
privately communicated
results 
on many occasions
\cite{BI95,HP91,Ne86,PD83,Pu85,WP87}. 
Thus,
approximate  analogs of 
 $++$
 were  explored 
subsequently
 by  Buss-Ignjatovic,
H\'{a}jek 
and
\v{S}vejdar in \cite{BI95,Ha7,Sv7},
as well as in Appendix A of 
our paper
\cite{ww1}.
Also, 
Pudl\'{a}k's initial 1985 article  \cite{Pu85} 
captured
the majority 
of $++$'s 
essence, chronologically before Solovay's observations.
Also,
Friedman did
related work
 in
\cite{Fr79a}.

\item[   $~~$b.$~~$  ]
Part of what makes  $++$ interesting is that 
\cite{ww1,ww5,wwapal,wwiqfs}
explored two methods for 
generalized arithmetics
to confirm their own consistency, whose
natural hybridizations  are precluded by $++$.
Specifically,  these results involve using
Example \ref{nex-2.3}'s 
self-referencing  ``I am consistent''
 axiom (from its
statement  $\oplus$ ). 
They will enable
several  Type-NS basis 
systems \cite{ww1,wwapal,wwiqfs}
to verify their   own consistency under 
a Hilbert-style deductive apparatus
\footnote{ The Example \ref{nex-2.1}'s 
provided
three examples of
  Hilbert-style 
deduction operators, called 
$\, d_E \,$, $\, d_H \,$ , 
and
 $\, d_M ~~$. It explained how these 
 deductive operators differ from a tableaux-style
deductive apparatus by containing a modus ponens rule.},
or alternatively allow 
 a Type-A 
 system \cite{ww93,ww1,ww5,ww6,ww14} to 
corroborate
 its  own 
self-consistency
under a more restricted semantic
tableaux style deductive apparatus.
Also, Willard \cite{ww2,ww7} observed how one could
refine $++$ with Adamowicz-Zbierski's
methodology \cite{Ad2,AZ1} to show  Type-M  systems
cannot recognize their own tableaux-style consistency.
\ennd
\end{remm}

\section{General Perspective}
\label{s3}
This section will explain how some seemingly 
minor
hair-thin
Boundary-Case exceptions to the Second Incompleteness Theorem
can be transformed into major chasms when one contemplates
the facts that 21-st century computers can perform arithmetic computations
 with more than a
billion times the speed
of the human mind (with a similar accompanying increase
in the lengths of the logical sentences being 
manipulated).
This distinction will raise questions about whether  21-st
century
engineering projects
 will ultimately be forced to encounter questions about
the Second Incompleteness Theorem, which were ignored when the
scientific world had
 first learned about G\"{o}del's  work 
during
 the
early
 1930's.

\subsection{Linguistic Notation}

\label{s3.1}

Our language
for formalizing exceptions to the Second Incompleteness theorem
will be called
 $~L^*~$.
It
 will include the
symbols $~C_0~$,  $~C_1~$
and  $~C_2~$ for representing the
integers of 0, 1 and 2. The language $~L^*~$ will
formalize other  non-negative integers 
using growth function primitives and these three
starting integers.

The only predicate symbols used by our language $~L^*~$ will
be the
equality and less-than-or-equal predicates,
denoted as $~$``$~=~$''$~$
and  $~$``$~ \leq~$''$~$.
Sometimes, we will informally also use the symbols
$~\geq~$, $~<~$ and $~>~$.

Define  $ F(a_1,a_2,...a_j) $ to be a {\bf NON-GROWTH FUNCTION}
iff for all values of $a_1,a_2,...a_j$, the
function  $F$ satisfies 
$~ F(a_1,a_2,...a_j) ~ \leq ~  Maximum(a_1,a_2,...a_j)~$. 
Our axiom systems will  employ a set of
eight non-growth functions, called the
{\bf GROUNDING FUNCTIONS}. They will
include:
\begin{enumerate}
\baselineskip = 1.2\normalbaselineskip
\item  Integer Subtraction where 
$ x-y$ is defined to equal zero when $x<y,~~$
\item  Integer Division where 
$~  \frac {x} {y} ~=~ x~$
when $~y=0~$, and it otherwise equals
$\, \lfloor \, \frac {x} {y} \, \rfloor \,$.
\item 
$Predecessor(x) ~=~\mbox{Max}(\,x-1\,,\,0\,),~~$
\item
$~Maximum(x,y),~~$
\item
$ Logarithm(x) ~=~ \lceil  Log_2 (x+1) \rceil,~~$
\item 
$Root(x,y)$ = $ \lfloor ~x^{1/y}~ \rfloor$ when $~y \geq 1~$,
and $~Root(x,0)~ =~x. $ 
\item
$Count(x,j)$ designating the number of ``1'' bits among
$~x$'s$~$ rightmost $~j~$ bits.
\item
$Bit(x,i)$
 designating the value of the integer $x$'s $i-$th rightmost
bit. (Note that
\newline
 $Bit(x,i)~=~Count(x,i)-Count(x,i-1)$.) 
\end{enumerate}

In addition to the above non-growth functions, our language $~L^*~$
will employ two growth functions. They will correspond to
Integer-Addition and Double$(x)=x+x$. We will use the term
{\bf U-Grounding Function} to
refer to a function that is one of the eight Grounding
Functions or one of the operations of Addition and Doubling.

\bigskip

\begin{comm}
\label{f3.1}
\label{vvc.1.}
\rm
 We do not
technically need both the operations of
Addition and Doubling in our U-Grounding language $~L^*~$.
However,  it much is easier to encode large integers
when we have access to both these function symbols.
For example, any integer $~N>1~$ can be encoded by a term of
length $O(~$Log$(N)~)~$, using only the
constant symbol for ``1'',  when both the addition and doubling
function symbols are present. For instance, below is 
our binary-like encoding for the number eleven.
\beq
1~+~\mbox{Double}(~1+~\mbox{Double}(~\mbox{Double}(~1)))
\enq
Henceforth, the symbol $\overx{N}$ will denote such a binary-like
encoding for the integer $~N~$. (In the degenerate case where
$~N=0~$,  $\overx{0}$ will simply be defined as being the constant
symbol ``$~C_0~$'', that represents zero's value.) 
\end{comm}

\bigskip

\begin{deff}
\label{f3.2}
\label{vvc.2.}
\rm
We will follow mostly
conventional logic notation when discussing
the
 U-Grounding functions.
Thus, a {\it term} is defined to be a constant symbol, a variable
symbol or a function symbol (followed by some input arguments,
which are
similarly defined terms).
If $t$ is a term then the quantifiers in
$~ \forall ~ v \leq t~~ \Psi (v)~$ and $~ \exists ~ v \leq t~~ \Psi (v)$
will be called {\it bounded quantifiers}. These two wffs will be
semantically equivalent to the
respective
 formulae of
$\forall v~ ~(~ v \leq t~ \Rightarrow  ~ \Psi (v))~$ and
$\exists v~ ~(~ v \leq t~ \wedge  ~ \Psi (v))$.
A formula 
$\Phi$ will be called
$\Delta_0^*$ 
iff all its quantifiers
are bounded. Thus a 
$\Delta_0^*$ formula is defined to be a wff that is
{\bf any combination}
of atomic formulae 
(using our ten U-Grounding functions and the equals and $~\leq~$
predicates)
combined with bounded quantifiers and with
the boolean operations of AND,OR, NOT and IMPLIES
in an arbitrary manner.
\end{deff}

\begin{deff}
\label{vvc.4.}
\label{f3.3}
\rm
For any integer $~i \geq~0~$, this paragraph will define
the notions of a 
$\Pi_i^*$ and  $\Sigma_i^*$ formula. Their definition has three parts,
and it is given below.
\bee
\item
Every $~\Delta_0^*~$ formula is defined to be also
both a $\Pi_0^*$ and  $\Sigma_0^*$ formula. 
\item
 A formula is called
$\Pi_{i+1}^*$ 
iff for some 
 $\Sigma_i^*$ formula 
 $ ~ \Phi(v_1,v_2,...v_n),~$ 
it can be
written in the form 
$\forall v_1 ~ \forall v_2 ~...~ \forall v_n ~ \Phi(v_1,v_2,...v_n)~.~~$
(Since this rule also applies when the number of quantifiers $~n~$
equals zero, it follows that every 
 $\Sigma_i^*$ formula is automatically by default also
$\Pi_{i+1}^*$.)
\item
Similarly,
a formula is called
$\Sigma_{i+1}^*$ 
iff for some 
 $\Pi_i^*$ formula
 $ ~ \Phi(v_1,v_2,...v_n),~$ 
it can be
written in the form  
$ \exists v_1 ~  \exists v_2 ~...~  \exists v_n ~ \Phi(v_1,v_2,...v_n)~.~~$
(Since this rule
 applies when the number of quantifiers $~n~$
equals zero, it follows that every 
 $\Pi_i^*$ formula is
also
$\Sigma_{i+1}^*$.)
\ene
\end{deff}

\begin{exx}
\label{f3.4}
\label{vvc.5.} $~$
\rm
Lines
\eq{e1} are \eq{e2} 
are examples of $\Pi_2^*$ sentences, and 
Lines
\eq{e3} are \eq{e4} 
are examples of  $\Pi_1^*$ sentences. Note that some   
 $\Pi_2^*$ sentences can be proven to be logically equivalent to
$\Pi_1^*$ sentences. Thus for example, the sentences in Lines
\eq{e1} and \eq{e4} are
equivalent to each other.

{\baselineskip = 0.7 \normalbaselineskip  \parskip   0pt
\begin{equation}
\label{e1}
\forall~ x~~~ \forall~ y~~~ \exists~ z~~~~~ ~z-x~=~y~ 
\end{equation}

\begin{equation}
\label{e2}
\forall~ x~~~ \forall~ y~~~ \exists~ z~~~~~
x>0~ \Rightarrow ~\frac{z}{x} ~=~ y
\end{equation}

\begin{equation}
\label{e3}
\forall~ x~~~ \forall~ y~~~ ~~~~~ x+y ~=~ y+x
\end{equation}

\begin{equation}
\label{e4}
\forall~ x~~~ \forall~ y~~~ \exists~ z \leq x+y~~~~~ ~z-x~=~y~ 
\end{equation}}
Note that while \el{e1}'s
$\Pi_2^*$ sentence 
 can be transformed into an
equivalent
$\Pi_1^*$ formula
(encoded by  \el{e4} ), there is no analogous
$\Pi_1^*$ sentence
which is
 equivalent  to  \el{e2}. (This is because Addition
{\it but not Multiplication}
 belongs to our
particular specified
 set of U-Grounding functions.)
\end{exx}

\bigskip

\begin{remm}
\label{f3.5}
\label{vvc.10}
\rm
Throughout all our papers, the symbols
Add$(x,y,z)$ and  Mult$(x,y,z)$ 
will
denote two
 $\Delta_0^*~$
 formulae that are satisfied precisely when
the respective  conditions of $~x+y=z~$  and  $~x*y=z~$
are true.
It turns out that we can define both these formulae
using only the non-growth functions of
integer subtraction and division.
Such
definitions
of  Add$(x,y,z)$ and  Mult$(x,y,z)$ are provided by Lines
\eq{adde} and \eq{multe} below
\beq
\label{adde}
z-x=y~~\wedge~~ z \geq y
\enq
\beq
\label{multe}
[~(x=0    \vee    y=0 ) \Rightarrow z=0~ ]~ ~\wedge ~~ 
[~(x \neq 0 \wedge y \neq 0~) ~ \Rightarrow ~
(~ \frac{z}{x}=y  ~\wedge \, ~  \frac{z-1}{x}<y~~)~]
\enq
In this context, an axiom system $~\alpha~$ will be said to  recognize
Addition {\bf ``as a total function''} iff it can prove
\beq
\label{addt}
\forall ~x~~\forall ~y~~\exists ~z~~~~ \mbox{Add}(x,y,z)
\enq
Likewise, we will say an axiom system $~\alpha~$ can recognize
Multiplication {\bf ``as a total function''} iff it can prove
\beq
\label{mult.t}
\forall ~x~~\forall ~y~~\exists ~z~~~~ \mbox{Mult}(x,y,z)
\enq
Also,
we will say an axiom system $~\alpha~$ can recognize
``Successor''
 as a total function iff it can prove:
\beq
\label{suct}
\forall ~x~~~\exists ~z~~~~ \mbox{Add}(x,1,z)
\enq

Some axiom systems $~\alpha~$ are unable to prove Multiplication is
a total function,
 but they can prove every true $\Delta_0^*$ sentence.
Other axiom systems 
are
unable to recognize any of Addition, Multiplication
or Successor as total functions,
but they can still  prove every true 
 $\Delta_0^*$ sentence. It will turn out these facts will be central
to understanding the generality and limitations of G\"{o}del's 
Second Incompleteness Theorem.
\end{remm}

\medskip

\begin{deff}
\label{f3.6}
\label{vvc.11} 
\rm
A sentence $~\Phi~$ will be said to be written
in {\bf  Prenex* Normal Form} iff for some $~i \geq 0~$, it can
be written as  a $~\Pi_i^*~$ or a $\Sigma_i^*$ sentence.
(It can be easily established that a predicate logic, using the language
$~L^*~,~$ can  show
 every 
 sentence
$~\Phi~$ can be mapped onto a Prenex* sentence $~\Phi^*~$ such that
 $~~\Phi ~ \Leftrightarrow ~ \Phi^*$.
Thus without any loss in generality, we may assume that all the proper axioms
within
a
basis system $\alpha$ 
can be
 encoded in a 
 Prenex* normalized form.)

\end{deff}


\subsection{Enriched Forms of Tableaux Deduction}

We will first employ our preceding language $~L^*~$ to review the definition
of ``Semantic Tableaux'' deduction in this section. We will then define two
minor variations of this construct,
 called the ``Rank-Zero'' and ``Rank-Zero-Plus''
enriched versions of
Tableaux deduction.


Our
definition of a
semantic tableaux proof
will be 
similar to its counterparts in
Fitting's and Smullyan's textbooks  \cite{Fit,Smul}. 
  Define a
{\bf $\Phi$-Focused Candidate Tree} for
the axiom system $~\alpha~$
to be a  tree structure 
whose root corresponds to
the sentence $~\neg \, \Phi~$
 and whose all  other nodes are
either 
formal
axioms of $~\alpha~$ or deductions from higher
nodes of the tree. 
 Let  the notation
``$~ \cal{A} ~ \Longrightarrow ~ \cal{B} ~$''  indicate 
 $~  \cal{B} ~$ is a valid deduction
when $~ \cal{A} ~$ is an ancestor of $~  \cal{B} ~$.
In this notation, the deduction
rules allowed in a candidate tree are:
\begin{enumerate}
\small
\baselineskip = 1.11 \normalbaselineskip
\item $~ \Upsilon \wedge \Gamma \, ~ \Longrightarrow ~ \, \Upsilon
~$ 
and 
$~ \Upsilon \wedge \Gamma \, ~ \Longrightarrow ~ \, \Gamma ~$ .
\item $~ \neg  \,\neg \, \Upsilon ~ \Longrightarrow ~ \Upsilon~$.  
Other Tableaux rules for
the ``$ \, \neg \,$'' symbol  are: 
$~\neg ( \Upsilon \vee \Gamma ) ~ \Longrightarrow ~ \neg \Upsilon
\wedge \neg \Gamma$,
$ \, \neg ( \Upsilon \rightarrow \Gamma )  \,  \Longrightarrow  \,   \Upsilon
\wedge \neg \Gamma \, $,
$ ~~~~\, \neg ( \Upsilon \wedge \Gamma )  \,  \Longrightarrow  \,  \neg
\Upsilon \vee \neg \Gamma \, $,
 $~ \,   \neg \, \exists v \, \Upsilon  (v)  \,  \Longrightarrow  \,  
\forall v   \neg \, \Upsilon  (v)  \, $ and
 $ ~\,   \neg \, \forall v \, \Upsilon  (v)  \,  \Longrightarrow  \,  
\exists v \,  \neg  \Upsilon  (v)$
\item A pair of sibling nodes $~ \Upsilon ~$ and $~ \Gamma ~$ is
allowed
when their ancestor is
$~\Upsilon \, \vee \, \Gamma~$.
\item A pair of sibling nodes $ \,  \neg \Upsilon  \, $ and $ \,  \Gamma  \, $ is
allowed
when their ancestor is
$ \, \Upsilon \, \rightarrow \, \Gamma$.
\item $~ \exists v \, \Upsilon  (v) ~ \Longrightarrow ~ \, \Upsilon(u) ~$
where $\,u \,$ is a newly introduced ``Parameter Symbol''. 
\item 
Our variation of Rule 5 for {\it bounded existential quantifiers}
of the form $~$``$~ \exists v \leq s ~$''$~$
is the identity:$~$
 $~ \exists v \leq s ~ \, \Upsilon  (v) ~~ \Longrightarrow ~ ~
u \leq s ~ \wedge~ \Upsilon(u) ~$
\item $\forall v \, \Upsilon  (v)  \,  \Longrightarrow    \, \Upsilon(t)  \, $
where $t$ denotes any U-Grounded 
term. These terms are defined to be
  parameter symbols,
constant symbols,
or  U-Grounding functions with  
recursively defined
inputs.
\item 
Our variation of Rule 7 
for
{\it bounded universal quantifiers} of the form
 ``$~ \forall v \leq s ~ $''
is the identity:$~$
 $~ \forall v \leq s ~ \, \Upsilon  (v) ~~~~ \Longrightarrow ~~~ ~
t \leq s ~ \rightarrow~ \Upsilon(t). ~$ 
\end{enumerate}
Define a particular leaf-to-root branch in a candidate
tree $~T~$ to be {\bf Closed} iff it contains both some sentence
$~ \Upsilon ~$ and its negation $~ \neg \, \Upsilon ~$.  
 A {\bf  Semantic
Tableaux} proof of $~\Phi~$  will then be defined to be
a candidate tree,
{\it all of whose  root-to-leaf branches are
closed,} such that the tree's root stores the sentence
$~ \neg \Phi~$ 
and where  all its other nodes are
either axioms of $~\alpha~$ or deductions from higher nodes.

\begin{deff}
\label{def3.7} 
\rm
Let $~Z~$ 
denote an arbitrary set of sentences in our language $~L^*$.
Recall that a node in a semantic tableaux 
deductive
proof from an axiom system
$~\alpha~$ is allowed to include any axiom of $~\alpha~$
as 
one of
its stored sentences.  Such a proof will be called a
{\bf Z-Enriched} proof if it
may
 also include any version of
\eq{xclm}'s
formalization
 of the ``Law of the Excluded Middle'' 
as a permissible logical axiom
for every
$~\Psi~\in~Z~$.
\beq
\label{xclm}
\Psi ~~\vee ~~ \neg ~\Psi
\enq
\end{deff}

It is well known that semantic tableaux proofs satisfy 
G\"{o}del's Completeness Theorem
 \cite{Fit,Smul}.
 This implies that the set
of theorems that 
are proven
 from 
an axiom system $~\alpha~$ via
a conventional (unenriched) version of semantic tableaux deduction
is the same as the set of theorems 
proven
from a Z-enriched version of this deductive
mechanism. 
Our main result in this section will show,
 however, that such proofs can have
their efficiency 
often
exponentially
shortened under
such enrichments, 
where \el{xclm}'s schema is treated as a
set of
 logical axioms (rather than as a collection of 
derived theorems).

\begin{deff}
\label{def-3.8}
\rm
Let $~\alpha~$ be an arbitrary set of proper axioms and $~D~$ denote
a deduction method.  
In this case, an
 arbitrary theorem $~\Phi~$ of $~\alpha~$ shall be
said to satisfy a  {\bf $\Psi$-Based Linear Constrained Cut Rule}
iff $~\Phi\,$'s 
shortest proof from  $~\alpha~$ (via $D$'s apparatus) is 
guaranteed to be
no longer than
proportional to the
 sum of the  lengths
of the  proofs of   $\Psi$
and   $\Psi~\Rightarrow ~\Phi$
from $~\alpha~$.
\end{deff}

\begin{exx}
\label{ex-3.9}
\rm
All Hilbert-style deduction methods (including Example \ref{nex-2.1}'s 
$d_E$, $d_H$, and $d_M$ Hilbert-style methodologies) 
employ 
 Linear Constrained Cut Rules for any arbitrary input sentence $\Psi$
(on account of the presence of their modus ponens rules).
This is the intuitive reason that
 Theorem $\, ++ \,$'s generalization of the Second Incompleteness
Theorem 
(from Remark \ref{nremm-2.5} )
applies to them.
We will soon see that
a
 similar generalization of the Second Incompleteness
Theorem applies to modifications of  semantic tableaux deduction
(when \el{xclm}'s invocation of the Law of the Excluded Middle
is
 available 
as a logical axiom 
{\it for every input sentence} $\Psi$).
\end{exx}

\begin{lemm}
\label{lem-3.10}
\rm
Let $~D_{\Psi}~$ 
denote an ``enriched'' deduction method, identical to
semantic tableaux deduction, except that 
\el{xclm}'s version of the Law of the Excluded Middle is available
as a logical axiom under  $~D_{\Psi}~$.
Then for an arbitrary theorem $\Phi$,
a $\Psi$-Based Linear Constrained Cut Rule
will be satisfied by the deduction method $~D_{\Psi}~$.  
\end{lemm}

{\bf A Brief Sketched Justification:}
The germane proof $~p~$ of $~\Phi~$ will follow the usual semantic
tableaux format by having its root
 store
the sentence
 $~\neg~\Phi~$.
The child of this root will
then
 consist of
\el{xclm}'s version of the Law of the Excluded Middle.
Also, the
two children
of this node
 will 
consist of the sentences of $~ \neg \Psi~$ 
and
$~\Psi~$. We omit the details,
but it is easy to then verify that:
\bed
\item[  a. ]
 one may insert a proof-subtree
below  $~ \neg \Psi~$  that is no longer than
linearly proportional to the length of 
 $~  \Psi\,$'s proof.
\item[  b. ]
 one may insert  a proof-subtree
below  $~  \Psi~$  that is no longer than
linearly proportional to the length of the proof of
$~\Psi ~ \Rightarrow ~\Phi~$.
\ennd
These constraints imply that the
 $\Psi$-Based Linear Constrained Cut Rule
will be satisfied.  $~\Box$.

\begin{deff}
\label{def-3.11}
\rm
There will be several types of ``Enrichments'' of the semantic
tableaux deduction method that we will 
examine
in the context of
Definitions \ref{def3.7}  and \ref{def-3.8}. These will include:
\bee
\item
 {\bf Infinitely Enriched} formalisms that allow
\el{xclm}'s variation of the ``Law of the Excluded Middle'' to 
become a logical axiom,
for any sentence $~\Psi~$ from $~L^*~$'s language. 
\item
 {\bf Rank-k Enriched formalisms} that allow
\el{xclm}'s variation of the ``Law of the Excluded Middle'' to 
be a logical axiom
when $\Psi$ is
 any  $\Pi_k^*$ or  $\Sigma_k^*$ sentence.
\item
 {\bf Rank-Zero Enriched formalisms} that allow
\el{xclm}'s variation of the ``Law of the Excluded Middle'' to 
be a logical axiom
when $\Psi$ is
 any  $\Delta_0^*$ sentence.   
\item
 {\bf Rank-Zero-Plus Enriched formalisms} that 
are a slightly stronger version of the Rank-Zero formalism
that take \eq{zplus} as a logical axiom for any
  $\Delta_0^*$ formula $\psi(x)$.
\beq
\label{zplus}
\forall ~x~~~ \psi(x) ~ \vee ~\neg ~ \psi(x)
\enq
\ene
\end{deff}

\begin{remm}
\label{rem-3.12}
\rm
Let $~\alpha~$ denote any axiom
basis
 that
includes the ten U-Grounding symbols.
Then if $~D~$ denotes the semantic tableaux deductive methodology
and if  all of $\alpha$'s axioms hold true in the standard model,
it will follow that
 \cite{ww5}'s IS$_D(\alpha)$  
formalism
will
be a self-justifying  system which  proves all $\alpha$'s
$\Pi_1^*$ theorems and additionally recognizes its own semantic tableaux consistency.
It turns out this result will also generalize when D corresponds to
either Item 3's Rank-Zero  enriched form of semantic tableaux deduction
or   Item 4's  Rank-Zero-Plus  enriched form.
 (E.g. these two systems will
be
also
 capable of
 recognizing their
own
 consistency under their enriched deduction
methods.)

 In contrast, one may easily 
apply Lemma \ref{lem-3.10} 
to show that
 the invariant ++
(appearing in Remark \ref{nremm-2.5})
will generalize to establish that Type-S axiom systems
cannot verify their own consistency under infinite enrichments of
the semantic tableaux
\footnote{This is because infinite enrichments of tableaux
make its deductive procedure resemble Hilbert deduction
because such enrichments uniformly obey
Definition \ref{def-3.8}'s 
 Linear Constrained Cut Rule.
In particular, these enrichments become inconsistent
if they employ an analog of
\cite{ww5}'s Group-3
{\it ``I am consistent''} axioms}.
%
%
Indeed, the methods from \cite{wwlogos} imply 
the
Second Incompleteness Theorem  also generalizes for the case of
Rank-2 
or higher
enriched formalisms
under Type-A arithmetics.
 Thus, there is only a meaningful
open  question about the
application of  Remark  \ref{rem-3.12}'s   
paradigm to
               Rank-1 Enriched systems.
\end{remm}

An attached appendix will review \cite{ww5}'s definition of the
 IS$_D(\alpha)$ axiom system, for the benefit of those readers  who
have not read 
\cite{ww5}.
(It will amplify upon the claims made in the previous
two paragraphs.)
Our recommendation is  that the
reader postpone examining this appendix until after the main sections
of the current paper are finished.
This is because Definition 
\ref{def-3.11}'s notion of ``Deductive Enrichment''
is 
quite
subtle, and the next 
two
sections
shall
 need to
first
 consider it 
in more detail.

\section{The Significance of Deductive Enrichment}
\label{sign}

A general rule of thumb in Proof Theory is that an axiom system is
typically extended 
in order to
 expand the class of theorems it can
prove.  Since semantic tableaux deduction satisfies G\"{o}del's
Completeness Theorem, the Definition \ref{def-3.11}'s four variations
of deductive 
enrichment
will, however,
 not 
technically
change the theorems that
 can be derived from
an initial base axiom system.  Instead, the function
of  deductive 
enrichment
 will be to
improve the
overall
 proof efficiency. (This is 
because an invoked version of
the Law of the Excluded Middle will shorten
proof lengths
when it is treated as a
system of
 {\it  logical axioms,} rather than as a set of 
{\it derived theorems.)
}

This issue was perhaps not so central in the early 20-th
century when G\"{o}del announced his initial 2-part
Incompleteness paradigm. At that time, the only available
medium of thought was
 the Human Brain, which  
performed arithmetic computations
at approximately
a billion times a slower
speed
 than 
that of the
typical 21st century   household computer.
Also, the potential lengths of logical sentences during the 1930's
was much shorter than
many potentially lengthy
 present-day computer-generated sentences.

Within the context of the  longer expressions
that computers 
can physically
 produce,
the task of separating true from false $\Delta_0^*$
sentences
will likely
  become
increasingly
 daunting (assuming $ P \neq NP$ ), even when  this task
is
 technically
decidable.
 Hence a Rank-Zero  Enrichment of a tableaux deductive
formalism is a useful 
instrument,
 with the improved efficiency of its
Rank-Zero Linear Constrained Cut Rule. 
Moreover, the results from our earlier paper \cite{ww5}
do trivially imply that our Rank-Zero and 
Rank-Zero-Plus extensions of Self-Justifying formalisms
are guaranteed to be consistent.
(This is because an
 application of 
either the Rank-Zero or Rank-Zero-Plus versions of
the Law of the Excluded Middle can be
formally encoded within a $\Pi_1^*$ format, and
the attached Appendix explains
 our results
from \cite{ww5} imply that the addition of
any set of
 logically valid
$\Pi_1^*$
sentences are
 compatible with IS$_D(\alpha)$
 retaining its internal
consistency.) In contrast, \cite{wwlogos} 
implies
 the same is not true
for Rank-2 and higher enrichments of tableaux deduction; e.g. 
these formalisms lose their self-justification property
when their 
 enrichments are incorporated.


\begin{remm}
\label{rem-4.1}
\rm
Many readers will be initially disappointed that Rank-2 and other
higher enrichments levels will be
 infeasible under self-justifying  semantic tableaux deductive systems.
This will mean that if $~\Psi~$ 
denotes \eq{tmult}'s 
declaration that multiplication is a total function
 then
neither
can it  be assumed to be
 true
by our self-justifying formalisms,
nor can 
the theorem  $~\Psi~\vee ~ \neg \Psi~$ be promoted
into becoming a {\it logical axiom} (under 
the
self-justifying methods of  IS$_D(\alpha)$
without producing an inconsistency).
\beq
\label{tmult}
\forall x ~\forall y ~\exists z ~~~Mult(x,y,z)~
\enq 
Nevertheless,
there is a method whereby our Rank-Zero Enriched
formalisms can partially formalize 
\el{tmult}'s meaning. 
Thus,  let $~\Psi_k~$ denote the $\Delta_0^*$ formula
(shown below)
indicating that its localized version of multiplication is a total function
among 
input
integers less than $~2^k~$.  
\beq
\label{kmult}
\forall x~ \leq~\mbox{Double}^k(2)~~~
 ~\forall y ~ \leq~\mbox{Double}^k(2)~~~
\exists z~ \leq~\mbox{Double}^{2k}(2)~~~~
 ~~~Mult(x,y,z)~
\enq 
Then it turns out that our Rank-Zero enriched versions of tableaux deduction
can  prove $~\Psi_k~$ 
as a theorem, as well as treat
 $~\Psi_k~\vee ~ \neg \Psi_k~$ 
  as
a logical axiom. We will call  $~\Psi_k~$  and the logical
axiom  $~\Psi_k~\vee ~ \neg \Psi_k~$
 the {\bf K-Localizations}
of the sentences 
 $~\Psi~$  and 
  $~\Psi~\vee ~ \neg \Psi~$.
\end{remm}

In many pragmatic applications, one does not
technically
 need 
$\Pi_j^*$ and $\Sigma_j^*$
theorems $\Phi$
(with $j \geq 1$): 
Instead,
it suffices
to prove a
 K-Localized
theorem $\Phi_k$,
that employs analogs 
of  \el{kmult}'s three specified bounded quantifiers,
for some sufficiently large fixed constant $~k~$.
In particular, a transition of higher sentences into 
$\Delta_0^*$ formulae is especially pragmatic 
for
21st century computers, whose speed and 
allowed byte-lengths can exceed by factors of
many millions their counterparts produced by human mind.

Within such a context, the
self-justifying
 capacities of even Rank-Zero Enriched
forms of Semantic tableaux deduction will be much more tempting
during the 21st century than they were at the time of 
G\"{o}del's 1931 discovery (when  computers 
were
 unavailable).
It is mainly for this reason that we suspect the modern world
should not fully dismiss the capacities of self-justifying logic
formalisms.

Moreover, we suspect that the futuristic civilization 
within our solar system,
 including
that
 on the planet
Earth, may have no choice but to rely upon 
Self-Justifying computerized logic systems.
This is because many scientists (including the late physicist
Stephen Hawking)
have
 speculated that if current trends continue, then
Global Warming will cause the planet Earth to become too hot
for mammals to survive on it, within one or two centuries.
In such a context, where  computers will not need the
Earth's 
 cooler temperatures
and/or Oxygen to survive
(e.g. see footnote 
\footnote{In particular, solar powered computers, 
physically
residing outside
the planet Earth,
 will
 need 
 neither
Oxygen as an energy supply,
nor have their operative functions
being
governed by the Earth's 
temperature
(when they lie
 physically
 outside the planet's
domain).
\smallskip
 } ),
Stephen Hawking
\cite{Ha17a,Ha17b}
has predicted that
computers
may become 
the
main
form of thinking agent 
during what will hopefully be 
only
a temporary
period
 of Global Warming.

%


Such computers will
need, presumably, 
 to rely upon some form of Self-Justifying
cognitive process to organize the motivations of their thought processes.
In particular, humans seem to  have relied upon some type of instinctive
appreciation of the coherence of their thought processes, as a 
prerequisite for motivating their cognitive process.
Our suspicion is that
computers will need to imitate this self-reinforcing process.

Our conjecture is,
thus, that the continuation of human civilization, 
within our Solar
System, may require computers taking 
temporary
control
of 
the larger part
 of
its destiny.
More specifically,
we suspect that a
computer-and-robot
technology 
shall
be  needed to
 reverse Global Warming and  enable
a  saved
sample of
frozen mammal
 embryos 
to subsequently populate
the planet Earth, again.  


Some readers may shutter at the thought that
planet Earth could become temporarily uninhabitable during a
perhaps thousand-year era of Global Over-Heating. This
difficulty,
however, may
 actually amount to only be a
temporary phenomena, 
if computers and
robots can restore Earth into a more hospitable environment after 
a period of several thousand years
(and also  frozen human embryos are saved).

More precisely in a context where life has existed on Earth for 
approximately 4 billion years, our perspective is that the danger
posed by Global Warming
{\it  would be  tragic
but temporary,}
if robots and computers can reverse
Global Warming after a period of a few thousand years.  In contrast,
the implications of Global Warming would be much more
severe,
if either it cannot be reversed or no frozen embryos
are saved
before Global Warming occurs.
 It is for this reason that we
suggest it is
 imperative
that
a fleet
of self-justifying computers,
along
with accompanying
 robots and
saved
 frozen
embryos, be assembled
as at least a partial response
 to a potentially
very
 serious Global Warming crisis.
(See footnote 
\footnote{ It is
possible
 that if computer and robotics
technologies
do
 advance {\it quickly enough} than a  tragic
 Global Over-Heating can be 
 completely avoided.
Moreover,
Steven Pinker discusses in \cite{Pi18}
some technologies that could possibly 
eschew
 global warming.
Our point is,
however, that it would be wise  to also investigate the potential
application
 of
self-justifying formalisms because their use may be necessary
in the future for Earth's mammals to survive and prosper.
The Corollary A.3, in Appendix A, thus illustrates a
potential use of self-justifying computerized systems that
possibly
 could
 be
urgent.}
for 
some
clarifications about the nature
and limits
of this proposal.)





\section{On the Motivation for Writing this Paper}

\gvs
\label{s5}


The author of this article has published several articles about
generalizations and boundary-case exceptions for G\'{o}del's
Second Incompleteness Theorem \cite{ww93}--\cite{wwiqfs}, including
six papers that have appeared in the JSL and APAL. 
Unfortunately, the author
has experienced both a mild stroke and a mild heart attack
during the summer of  2016.  

These events did not
prevent the
author from continuing his teaching during 2017 and 2018.
They did, however, cause a change
in
 the
specific
 goals of our research program.
The
 current article 
was intended,
mostly, to
 encourage 
others to
join  this
futuristic
 research project.
It has, thus,  focused on explaining why
this topic warrants
further investigation.

This is subtle because our
IS$_D(\alpha)$ axiom
system
 (formally defined in the Appendix)
 has three disadvantages when $D$ denotes semantic
tableaux deduction. These drawbacks are that:
\bee
\item
IS$_D(\alpha)$ is
 an unorthodox
``Type-A'' axiom system, which recognizes
addition {\it but not also multiplication} as a  total function.
\item
IS$_D(\alpha)$ has employed a semantic tableaux deductive apparatus,
that
is 
substantially less efficient
 than a more conventional Hilbert-styled
deductive apparatus.
\item
IS$_D(\alpha)$
is able to recognize its consistency
only by employing a version 
of Example \ref{nex-2.3}'s 
self-referencing  ``I am consistent''
 axiom.
\ene

Our reply to Item 1 is that while 
IS$_D(\alpha)$ is
unconventional, it is {\it not quite as weak} as it may
first
appear. This is partly because  IS$_D(\alpha)$ 
treats multiplication as a 3-way predicate
Mult$(x,y,z)$, formalized by \el{multe},
rather than as a total function.  Moreover, 
 IS$_D(\alpha)$ can be  easily arranged to prove all
the $\Pi_1$ theorems of Peano Arithmetic (except for
minor changes in notation) when its input axiom system
$~\alpha~$ is made to correspond to 
the trivial extension of
Peano Arithmetic   that includes our ten U-Grounding functions
symbols.
While these adjustments may not be ideally Utopian,
they are
 sufficient to reply to
the main difficulties raised by
Item 1.

Our reply to Item 3 is also easy because the goal of 
IS$_D(\alpha)$ 
{\it IS NOT TO}
 prove its own consistency
{\it under 
the most challenging
 definition of a proof.}
Instead,
 it is  to
find an axiom system that is comfortable with a {\it built-in}
internal assumption about its own consistency
(via its 
physically
 built-in 
  ``I am consistent'' axiom).
Most conventional axiom systems do not
own extensions of themselves that
 {\it can  support even this task}
 (due to various extensions of the Second Incompleteness
Theorem that preclude the consistency of such extensions).
Fortunately, 
IS$_D(\alpha)$ is able to support this precious
self-justification property
when $D$ corresponds to either semantic tableaux deduction,
or one of its refined
Rank-Zero or Rank-Zero-Plus enrichments.

Many readers will probably  be
especially
 leery of
 Item 2 
because conventional Hilbert-style deductive proofs
satisfy Definition \ref{def-3.8}'s
Linear Constrained Cut Rule, while Semantic Tableaux
Deduction does not
obey this property.
In particular,
our IS$_D(\alpha)$ formalisms
are required to treat the Law of the Excluded Middle
as a
set of derived theorems (rather than as a
{\it formally
 stronger}
built-in
 schema
 of
logical axioms).
Such a constraint will
make our proofs
substantially
less efficient
 than a more ideal
paradigm centered around 
 Definition \ref{def-3.8}'s
Linear Constrained Cut Rule.

  We can,
fortunately,
partially reply to 
this
daunting
 challenge
 because 
 Definition \ref{def-3.11} indicated
that there were
 four
methods for enriching
semantic tableaux deductive machineries.
Two of these four methods (i.e. the  Rank-Zero and
 Rank-Zero-Plus enrichments)
are
formally
compatible
with 
the existence of
self-justifying extensions of their semantic tableaux
deductive machineries.
Thus at least the 
  Rank-Zero and
 Rank-Zero-Plus
versions of the Law of the Excluded Middle may
be incorporated into self-justifying 
 semantic tableaux formalisms.

Naturally, it would be better if
 Remark \ref{rem-3.12} indicated
enrichment methodologies of Rank-2 and higher were also
compatible with self justification.
Unfortunately, however, excessive enrichments of the 
semantic tableaux deductive systems preclude self-justifying 
systems from existing (in a manner
roughly
 analogous to how excessive
enrichments of Uranium Ore can cause  nuclear reactors
to become 
dangerously
unstable). Thus, we 
are currently
confined to study of Rank-Zero
and  Rank-Zero-Plus
 enrichments of semantic tableaux, 
in a context where Rank-2 enrichments are infeasible,
and  Rank-1 enrichments are a
remaining
 open question. (See \cite{wwlogos}
for a discussion of  Rank-2 enrichments.)


While one should not ignore the fact
{\it that only}
  Rank-Zero and 
Rank-Zero-Plus enrichments  of semantic tableaux
deduction are known to be compatible with self justification,
it should be
noted
 even such Rank-Zero
enrichments are interesting.  This is because
one can philosophically argue that a logical sentence
loses its purely sensuous quality when it employs
unbounded quantifiers.  Thus even $\Pi_1^*$ and
$\Sigma_1^*$ sentences lie slightly above the 
{\it ``touch-and-feel level''} on account of their use
of unbounded quantifiers.   In other words,
our available ability to muster self-justifying 
 Rank-Zero and 
Rank-Zero-Plus enrichments  of semantic tableaux
deduction is significant because this level of enrichment
is broad enough to include 
the
crucial 
{\it ``touch-and-feel''} sentences of the language $L^*$.

The preceding point is important because it allows us
to summarize both the strengths and weaknesses of 
G\"{o}del's 1931 Second Incompleteness Theorem.
Thus, the traditional literature has been
certainly correct
in viewing 
G\"{o}del's discovery as a seminal result, when
our boundary-case exceptions to it
have persisted 
at only 
the  Rank-Zero and 
Rank-Zero-Plus enrichment levels. On the other hand,
the existence of such enrichments show that
G\"{o}del and Hilbert were 
partially
 correct when their
statements $*$ and $**$ foresaw that some types
of exceptions to the Second Incompleteness Effect
would persist (e.g. see \textsection  \ref{nnn1} ).
Moreover,
 the study of 
how to efficiently process
$\Delta_0^*$ sentences is 
important both because of their
{\it ``touch-and-feel''} property
 and
because it will be 
challenging to 
 process  
these sentences
efficiently,
 assuming that 
$P \neq NP$. 

Moreover, we remind the reader that 
Stephen Hawking
 and other scientists 
have expressed
 concern that
Global Warming could possibly lead to,
at least,  a temporary era
where
digital computers 
replace human brains as the
primary economically efficient mechanism for generating
thoughts \cite{Ha17a,Ha17b}.
In a context where computers can generate thoughts
more than a billion times faster than the human mind,
the preceding chapter
 suggested that
a self-justifying formalism, using
only Rank-Zero and 
Rank-Zero-Plus enrichments 
of semantic tableaux,
could help  computers 
function
 more efficiently.

In particular, our 
fervent
hope is that humans will  continue
to be 
 important
and central
in the future. It is likely
that
computerized simulations of the human thought process
will be also   important,
even if Global Warming occurs
as,
hopefully, a very
temporary phenomena. 


The writing of this section has been 
certainly
a
 painful
task
because it conveys
a message that is an
awkward
 mixture of  insight, hope and 
humbling realization. Thus, the 
fundamental
 insight,
appearing herein,
 is that some types of
borderline exceptions to the Second Incompleteness Theorem 
are noteworthy because their
logical formalisms can retain at least some type of
partial internalized
appreciation of 
 their
own
 consistency. The then-accompanying hope
is that computers can use this knowledge to help
life on Earth survive Global Warming and other futuristic
challenges.  And finally our
 humbling,  albeit
partially
 optimistic, conclusion
is that 
a prosperous life can continue on planet Earth
if (?) and
when 
 humans share with
computers
a
joint  control over our future destiny.
(We wish we could offer here a more positive 
 conclusion, but this section
has 
offered
at least
a partially optimistic message).

\section{Extending Our Main Thesis}


\label{s6}

Each of \textsection    \ref{nnn1}-\ref{s5}
 had appeared in an
 earlier  July 2018 version of our Cornell Archives technical 
report.  They had
subsequently  undergone only
very minor revisions editorial revisions.
 In contrast, this additional chapter is new.
Its goal will be to show that self-justifying arithmetics
will have quite significant implications,
{\it even if Global Warming does not occur.}

\subsection{Background Setting}  


\label{s61}

Our prior research during the last 25 years had received a positive
reception. It thus generated six sole-authored papers
in the JSL and APAL.

At the same time, the logic community had also treated our
work with
certainly some
guarded
 caution. This is because we documented only the existence of
 boundary case exceptions to the Second
Incompleteness Theorem.
More specifically,
there can be no disputing
that  
the Second Incompleteness Theorem has
 shown that at least roughly 97 \%
of the initial goals of Hilbert's Consistency Program
were simply infeasible.

The significance of the prior sections of this report is
that they have shown how a nontrivial final 3  \%
of the implications of the Second Incompleteness Theorem
should be further investigated.
This is because one cannot
doubt that Global Warming is a 
certainly
potentially major threat.
We thus observed how self-justifying arithmetics can
represent,
at least, a {\it fractional} 
\footnote{Obviously, many other technologies  will be needed,
along with  Symbolic Logic,
to alleviate Global Warming. They
will include 
a more sophisticated understanding of
 robotics, artificial intelligence,
the cognitive sciences, 
the ability to engage in
space travel and to better engineer  biological systems.
Our point, here, is only that a more sophisticated understanding
of Self-Justifying logics 
will be
 one part of a much broader solution
to the threat of Global Warming.  }
 part of a solution 
to the Global Warming dilemma.


Most scientists suspect that  Global Warming is a real threat,
but we shall not 
include
 such an assumption within this chapter.
Instead, we will focus on the famous interstellar paradox that was
posed by Enrico Fermi during a 1950 luncheon at the
Los Alamos National Laboratories \cite{firmiwik}.
We will
 show that
self-justifying arithmetics are equally important within
the context of the 
perplexing
Fermi Paradox.   

\subsection{The Fermi Paradox}


\label{s62}

The Fermi Paradox
\cite{firmiwik} 
has been nicely summarized in Wikipedia
as
highlighting
 the essential contrast between arguments showing the scale and
probability of the existence of intelligent 
life being common in the universe versus the total lack
of evidence of life as arising anywhere outside   Earth.


Many logicians would probably prefer Fermi's open question
to be called a ``Mystery'' 
or an ``Unsettled Question'' 
because it does not quite amount to being a ``Paradox'',
as  logicians
commonly use this term.

In particular, Wikipedia estimates it would take 
a civilization
no  more than roughly  ``5-50''' million years
to colonize 
a quite substantial number of
Earth-like
planets
in the Milky Way (once 
and if (?)
life reaches 
planet Earth's 21-st century
stage of civilization
and continues to survive).  This is because even if life
takes approximately $~ R~$ years to spread from 
one star system X to another system,
called say  ``Y'',
then it should be able to multiply throughout  the Milky Way's
no more than its
approximately  $4 \, \cdot \, 10^{11}$ possible
life supporting star systems
in roughly $~R~\cdot~ \mbox{Log}_210^{12}~$ years. (E.g.
there should be needed
no more than
 $~40~$ cycles of doublings of the population
of civilized stars until the majority of
such  stars
become life supporting.) 

In other words, if $R$ was even the relatively large number of $10^5$ years,
then only $4 \, \cdot 10^6$ years would be sufficient for civilization
to spread
to 
a great number
Earth-like
planets in the Milky Way
(when only 40 doublings
are
 sufficient to colonize
 the entire
 galaxy's
collection of life-supporting stars under a reasonable approximation
of the exponential growth associated with what should perhaps be called
the {\it ``Fermi Exponential Fertilization  Paradigm''}). 
On the other hand, a 
slightly larger
 $4 \, \cdot 10^7$ years will be
needed
 if, say,  $R=10^6$ years. 

Fermi asked the question about why such a colonization has not already
occurred in the Milky Way when
many
 stars much older than the Sun are known
to exist, and 5 or 50 million years is a relatively small amount
of time (compared to the age of the Milky Way). We will not delve into
this question here because it is 
complicated.
(It includes the possibility that all the predecessors to Earth's
civilization
have been unable to preclude their own self-destruction. It also allows
 for
 the possibility that some enlightened civilization will  ultimately 
find a more amicable outcome.)

Our point is that this topic  will  intersect,
eventually,  with one's
interest about the long-range implications of the Second Incompleteness
Theorem.  This is  because the 
time delay
          between communications among   distant
stars will,
presumably,
tempt
 an evolving civilization
to ultimately use a common computerized algorithm
(rather than
to
 deploy
time-delayed 
 radio signals as   
a medium to communicate a common logical formalism
among  distant star  systems). 
 
Thus while the Second Incompleteness Theorem's prohibition
against a formalism
owning an ideally robust appreciation of its own consistency
is a very serious obstacle, some type of possibly
fragmented understanding of 
one's self consistency
 will 
probably
be necessary, if a civilization 
within our galaxy can survive the challenges that time
thrusts upon it
(irregardless of whether or not Global Warming shall occur).

In other words, this chapter has reinforced the theme of the
prior section, {\it without assuming the 
existence of a  Global Warming  threat}.  
Thus while the Second Incompleteness Theorem and
its various generalizations
do
 show that at least 97 \%
of the
initial
 goals of Hilbert's  Consistency Program
were infeasible, some  diluted 
version of Hilbert's  program
will likely  be necessary
 for the long-range
futuristic
 survival of
an on-going civilization.  

\section{Overall Perspective}

\label{s7}

It should be noted that we have described  an alternate version
of self-justifying arithmetics in \cite{wwiqfs}
that are Type-NS systems that use a Hilbert-style
deductive apparatus, instead of
a
 semantic tableaux deductive apparatus.
These systems will be required to replace the assumption that addition
is a total function with 
\cite{wwiqfs}'s alternate
{\it weaker} 
 growth primitive, called the
$~\theta~$ operator,
so that its systems may remain consistent.
 We have not discussed \cite{wwiqfs}'s results
here
 because 
these logics have
  relied upon a conjecture  we consider
 likely, although it 
currently
is an unproven
 hypothesis,

We are uncertain which of our main results in this article or
in \cite{wwiqfs} are preferable.  Our overall perspective is
that the  threat of global warming, combined with the additional
long-term
issues associated with the Fermi Paradox, will make it inevitable
for an
evolving cognitive
 community to contemplate with increased seriousness
the implications of some type of partial evasions of the Second
Incompleteness Theorem. 

Especially in the context where I have experienced
a mild  stroke and
a mild heart attack during 2016,  I am  eager to
encourage other researchers join 
us in
this on-going
futuristic
 research project.

As   has been  already noted,
 the writing of this article
has been a quite
 painful
experience
 because it is disconcerting to consider
a future where humans may be required to share with
computers a joint control over our collective destiny.
Yet, there appears to be no viable long-term alternative for
the human race's survival. The good news is that 
a wisely controlled human-computer interaction can
preserve at least the main spirit of what humans are trying
to
accomplish

{\bf Additional Remark:}
This paper has been written as an ``extended abstract'',
 rather than as a fully detailed exposition,
primarily because many of our other papers
have discussed
related topics in much
greater  detail.
The writing style of an   informal
 extended abstract
 also seemed preferable because our chief goal
has been
 to stimulate a 
larger 
community
 of researchers
to join us in an ambitious  project,
centered around
futuristic
 self-justifying logic systems.

\section{About the Author and the Broader Impact of 
this Research}

In addition to publishing several papers about symbolic logic,
Dan Willard has been also active in several other areas of
research.  This other activity has included:
\bed
\item[  A. ] 
Co-Authoring with Robert Trivers \cite{TW73,wwtriv}
 a hypothesis that Darwinian
evolution exercises partial control over the sex determination
of offspring. 
 Google Scholar has  recorded
 3,648 citations 
to the TW article
\cite{TW73}, as of November 6, 2018.
The significance of the TW Hypothesis
 was also 
mentioned
          twice in the New York Times
on February 17, 1981 and on  January 21, 2017 in
articles entitled
{\it ``Species Survival Linked to Lopsided Sex Ratios''}
and {\it ``Does Breast Milk Have A Sex Bias''}.
(Moreover, Trivers 
has cleared up a possibly
awkward
statement, he
has
made,
 when page 113 
in
 his
     recent
 book
\cite{Tr2002}
does state
 the seminal idea behind the Trivers-Willard
Hypothesis  was 
fundamentally due
  \footnote{ The exact quote appearing on Page 113 
in
Bob Trivers's book
\cite{Tr2002} 
is the
frank
statement
 that if
he ``wanted to divide up credit''
exactly
 between
Trivers and Willard, 
the closing sentence of
\cite{TW73}
 should have
declared that
``The  idea was Dan's; the exposition Bob's'',
rather than
the passage that was written.} 
 to
Dan Willard's initial observations.)
\item[  B. ] 
Co-Authoring with Michael Fredman
 the
 articles  \cite{FW93,FW94}
about Sorting and Searching.
These papers
 were the chronologically first
among six items mentioned in the ``Mathematics and Computer Science Section''
of the {\it 1991 Annual Report of the National Science Foundation}.
(The prominence of Fusion Trees  was also signaled by  the
facts that they are the topic of Lecture 12 in MIT's 
Advanced Data Structures course, taught by its famous
MacArthur-award
winning professor Eric Demaine
(disseminated over YouTube \cite{YouDemaine}),
and also the topic of 
 the Lecture 2 in Harvard's
course Compsci 224 (recorded  on YouTube by
Harvard's Prof. Jelani Nelson \cite{YouNelson}.)
\item[  C. ] 
Developing the Y-Fast Trie
data structure
 \cite{ww83},
the Down Pointer Method \cite{ww85}
 and many
other results germane to Computational Geometry,
 Advanced Data Structures and Relational Database
Optimization Theory \cite{wcp,wc0,wc1,wc2,wc3}.
\ennd
The reason we have mentioned Items A-C here is
that we are aware that Willard's prediction
that Self-Justifying axiom systems could possibly  alleviate
the Global Warming crisis will
likely
 be seen as controversial.
This is because the human personality, which is partially
a product of Darwinian Evolution, 
is biased and
 uncomfortable
with the idea of computers sharing with humans
a knowledge of one's own
self-consistency.   
(Such knowledge implies an
uncomfortably
 awkward
shared 
 control
over our  joint
 future destinies, between computers and
mankind.)

Our hope is that the momentum established from our previous
research will 
add credibility to this paper's
awkward
 but 
quite
realistic
 predictions.
We hope this will 
 make a larger audience recognize
the  provocative  need for  artificially intelligent computers
to share
with humans
a temporary control over the planet Earth's
destiny until it
recovers from a global warming crisis.
We suspect  this
will be
 necessary for planet Earth
to become
 more hospitable for human
civilization's  vibrant  and
continuing
 survival.

Thus, the purpose of Items A-C was not to boast or
brag. Instead, it was to tempt 
our
 readers to
realize that a  new step  may be needed for 
civilization and humans to continue to 
peacefully thrive
 on Earth,
under the shadow of a Global Warming's pending threat
\footnote{ 
 This step will  require, of course, much more than the
evolution of self-justifying arithmetics. It will also require
advances in robotics, space travel, biology,
the cognitive sciences, and
AI-like
heuristic optimizing methods. 
We suspect a conjecture we made about  \cite{wwiqfs}'s
special $\theta$ function will turn out to be  also 
useful in this endevour. }.
Thus in a context where the majority of readers are
likely to be
understandably concerned about a World where computers might share 
equal power with humans over
 our
collective
 destiny, we hope   the prior success
of  our projects, A and B,
together with the
alarming
 shadow of a
likely
 upcoming 
Global Warming threat and
the prospects of a troubling
Fermi Paradox, 
will
help  our readers
to gain
 a renewed optimism
and  hope
 about what 
shall be, possibly, the
final
 necessary step
 for the sake
of 
the human spirit's 
continuing 
and thriving
survival.

In essence the classic interpretation, where the
Second Incompleteness Theorem has been seen as refuting
roughly 95-98 \% of the
initial
 objectives
of Hilbert's Consistency Program,
is 
basicly
 correct ---
{\it except for the fact} that it has
  ignored a 
{\it quite crucial} point
about how
21st and 22nd century logics shall 
have 
very
different goals than
those that
were
 envisioned, initially,
during the first half of the 20th century.

\bigskip

{\bf Acknowledgment:} I thank Seth Chaiken for several comments about how
to improve the
style of
 presentation.

\newpage
\section*{Appendix}




Our article \cite{ww5} had introduced
 both
 the definitions of 
the IS$_D(\alpha)$ axiom system and of a formalism
owning a ``Level$(J)$ Appreciation'' of its own consistency.  This appendix will
review the definitions of 
 these 
two
concepts, as well as
 explain how they are 
related to enriched versions of semantic tableaux deduction.
This appendix should be read only 
%
after Sections 1-\ref{s5}  are completed

\bigskip

During our discussion, 
 $L^*$ will denote the language defined
in \textsection \ref{s3.1},
and
 $~\beta~$ will denote a 
basis axiom system
whose
 deductive apparatus is denoted as 
D.
 In this context,
the
basis
 system
 $\beta$ will be said to  own: 
\bee
\item
a {\bf Level(n) Appreciation} of 
its own consistency under $D$'s deductive apparatus
iff $\beta$ can prove 
that there exists no $~\Pi_n^{* }$ sentence
$~\Upsilon~$ such that
 $~( \beta ,D )~$
supports  simultaneous proofs
of
both 
$~\Upsilon~$ and $~\neg \, \Upsilon .$ 
\item
a {\bf Level(0-) Appreciation} of 
its own consistency under $D$'s deductive apparatus
iff it can prove there 
exists no proof of 0=1 from the
 $~( \beta ,D )~$ formalism
\ene
 All these definitions of consistency,
from Level(0-) up to Level(n),
  are equivalent to each other
under strong enough models of Arithmetic.
However, many weak axiom systems do not possess the mathematical
strength to                    recognize
 their equivalence.

In particular, an axiom system
$\beta$
 owning a
Level(1) appreciation
 of its own consistency
is much stronger than 
such a system possessing a 
Level(0-) appreciation of its own consistency.
This is
because Level(1) systems can use a proof of
a $\Pi_1^*$ theorem
$~\Upsilon~$ 
 to gain the practical knowledge
that no proof of  $~\neg \, \Upsilon $  exists.

In a context where $\alpha$ is essentially any axiom system
that employs $~L^*\,$'s language and where 
$~D~$ denotes any deductive apparatus,
our  axiom system IS$_D( \alpha)$ 
in
 \cite{ww5}
was designed
 to achieve two specific goals. 
These were:
\bee
\item
To prove all true $\Delta_0^*$ sentences, as well
as to
prove
 all the $\Pi_1^*$
theorems that are implied by the axiom basis $~\alpha~$.
In particular, the formal definition of 
 IS$_D( \alpha)$ system,
 in 
 \textsection 3 of \cite{ww5},
accomplished the first task 
trivially
via its
``Group-Zero'' and ``Group-1 schema''. It performed
the second task because its 
 ``Group-2 scheme''
employed
\el{group2}'s 
particular
generic structure for each
$\Pi_1^* $ sentence  $~\Phi~$
(e.g. see Footnote \footnote{The symbol
$~\lceil \, \Phi \,\rceil~$ in $\eq{group2}$
denotes $~ \Phi ~$'s G\"{o}del number,
and the symbol $~ \mbox{Prf}~_\beta(~\lceil \, \Phi \,\rceil~,~p~)~$ will
designate a $\Delta_0^*$ 
formula
stating that $~p~$ is
a proof of $\Phi$.} ).
\begin{equation}
\forall ~p~~~\{~ \mbox{Prf}~_\beta(~\lceil \, \Phi \,\rceil~,~p~)~~~
\rightarrow ~~~ \Phi~~~\}
\label{group2}
\end{equation}
\item
To
``formally recognize'' its own Level-1
consistency under $D$'s deductive apparatus.
This 
was
 accomplished 
 by having
 IS$_D( \alpha)$'s ``Group-3'' axiom 
formalize
a $\Pi_1^*$
 sentence that 
amounts to 
\eq{group3}'s 
 statement. 
(The symbol 
``Pair$(x,y)$''
in \el{group3} is
 a $\Delta_0^*$ formula
indicating that  $~x~$ is the G\"{o}del number of a
 $\Pi_1^*$ sentence and that
 $~y~$ represents
 $~x \,$'s
``mechanically \footnote{In particular, if $~x~$ is the G\"{o}del
number of a $\Pi_1$ sentence $~\Phi~$ then  
its ``mechanically formalized negation''
 $~y~$ is 
the G\"{o}del
number for ``$~\neg~ \Phi~$''.
} formalized negation''.
Also,
$~\mbox{Prf}~_{\mbox{IS}_D(\alpha)}(a,b)~$
in \eq{group3} 
 denotes
a $\Delta_0^*$ formula 
indicating that $~b~$ is  a
proof, using the deduction method $~D~~,~$ of the theorem $~a~$ from
the axiom system IS$_D(\alpha)$.)
The nice aspect of \eq{group3} is that 
 IS$_D(\alpha)$
 can unambiguously interpret
the meaning of its three $\Delta_0^*$ 
formulae because its Group-Zero
and Group-1 schema allow it to correctly decipher \el{group3}'s
 statement.
\begin{equation}
\forall  ~x~\forall  ~y~\forall  ~p~\forall  ~q~~~~ \neg ~~
[~~ \mbox{Pair}(x,y)~ \wedge ~ 
~\mbox{Prf}~_{\mbox{IS}_D(\alpha)}(x,p)~
\wedge ~ ~~\mbox{Prf}~_{\mbox{IS}_D(\alpha)}(y,q)~ ]
\label{group3}
\end{equation}
\ene

\smallskip

{\bf Remark A.1}
It is unnecessary to provide a formal
description of the axiom system 
 IS$_D( \alpha)$ 
here because 
 \textsection 3 of \cite{ww5}
already
explained how it (and especially 
\el{group3} ) can be exactly encoded.
There are two
clarifications relevant to
 IS$_D( \alpha)$'s 
definition that should, however,
be 
mentioned.
\bed
\item[   a. ]
Our  IS$_D( \alpha)$ axiom system is a well defined entity for any
deductive apparatus $~D~$ and 
for any
axiom basis system
$~\alpha~$ 
(that have  recursively enumerable
formal
 structures).
Since \el{group3} 
is a self-referencing sentence,
one needs some meticulous  care,
however,
 to assure that the Fixed Point Theorem
can encode a 
 $\Pi_1^* $ sentence that is equivalent to  
\el{group3}'s statement.
(We do not discuss this topic here because it
was
discussed in 
 adequate detail in \cite{ww5}. )
\item[   b. ]
Although the  IS$_D( \alpha)$
 axiom system is well-defined entity
for all inputs $D$ and $~\alpha~$,
this fact
{\it does not} also guarantee that
 IS$_D( \alpha)$ is consistent.  Indeed, the Second Incompleteness
 Theorem implies 
 IS$_D( \alpha)$  is inconsistent for most
inputs  $D$ and $~\alpha~$. This issue was previously
visited by us in Example \ref{nex-2.3}.
It 
had
emphasized
that
\newline
{\it ``I am consistent''}
 axioms, similar to \el{group3},
can be easily encoded (via the Fixed Point theorem).
However, these axiom sentences
are {\it typically useless}, on account of the
inconsistencies that they
usually
 produce.  
\ennd


Our main result in \cite{ww5} 
is related to Item (b)'s 
ironical
paradigm.
In particular,  \cite{ww5} established 
that 
Remark A.1.b's
 paradigm can be
evaded when $D$ corresponds to
semantic tableaux deduction.
Thus,
 \cite{wwis,ww5}
 established that
 the
following partial boundary case exception to the
Second Incompleteness Theorem does arise:

\medskip

{\bf Theorem A.2.}  \cite{wwis,ww5}
 Let $~D_S$ denote the semantic tableaux deductive
apparatus, and $~\alpha~$ denote any axiom system all of whose
$\Pi_1^*$
 theorems 
are true sentences under the standard model
using the language $~L^*~$. Then   
 IS$_{D_S}( \alpha)$ will prove 
all $\alpha$'s
$\Pi_1^*$ theorems, and it will
also be consistent. 

\medskip

We remind the reader that Theorem A.2
 is  significant
because
 the Second Incompleteness Theorems implies
that most systems, formally verifying their own consistency,
actually
 fail to be consistent. 
 Theorem A.2
 is germane to 
 the current article 
because it implies the
following corollary:

\medskip

{\bf Corollary A.3} Let $~D_E$ denote the
 Rank-Zero-Plus
 enrichment of the
 semantic tableaux apparatus.
Also, let $~\alpha~$ 
again
denote an axiom system all of whose
$\Pi_1^*$ theorems 
are true sentences under the standard model.
 Then   
 IS$_{D_E}( \alpha)$ will also
 be consistent. 

\medskip

Corollary A.3 
is an easy consequence of Theorem A.2.  We will
now present a brief sketch of
its proof.

\medskip

{\bf Proof Outline:} Let $\alpha^*$
 denote the extension
of the basis system $\alpha$ that includes
one instance of axiom \eq{axor}
for 
every $\Delta_0^*$ formula $~\Psi(x) ~$.
\beq
\label{axor}
\forall ~x~~~~~ \Psi(x) ~~ \vee ~~ \neg ~ \Psi(x)
\enq
Theorem A.2 implies  that all the theorems of  $\alpha^*$
hold true under the standard model
simply
 because 
\eq{axor} holds true under the Standard Model,
as do
also
 all the 
particular
$\Pi_1^*$
 theorems of $\alpha$.
Hence, Theorem A.2 implies that
 IS$_{D_S}( \alpha^*)$ is consistent.  
This implies that 
  IS$_{D_E}( \alpha)$ must also be consistent
(since the 
 IS$_{D_S}( \alpha^*)$ and 
  IS$_{D_E}( \alpha)$ formalisms are 
 essentially
identical to each other, except for  minor changes in
notation). $~~~\Box$

\medskip

{\bf Remark A.4.}
The simplicity of Corollary A.3's proof
may tempt one to 
partially
overlook its
significance.
This corollary is significant because  
one of the main themes of our article has been
that a deductive apparatus does not capture 
the core intentions of
most  logics, unless it contains
some form of Definition \ref{def-3.8}'s 
 Linear Constrained Cut Rule.
The significance of Corollary A.3
is
that
 it shows
that  the
 Rank-Zero-Plus variant  of
 Definition \ref{def-3.8}'s 
 linear cut rule is
actually formally supported 
by the  self-justifying
IS$_{D_E}( \alpha)$
 axiom system.

{\bf Remark A.5.}
We again remind the reader that the formalism 
 IS$_{D_S}( \alpha)$ and
 IS$_{D_E}( \alpha)$ 
in Propositions A.2 and A.3 are capable of proving all
Peano Arithmetic's $\Pi_1^*$ theorems when
$~\alpha~$ designates
the trivial extension of
Peano Arithmetic   which includes the ten
 U-Grounding functions
symbols.
(These
self-justifying
 formalisms can thus appreciate 
a non-trivial part of 
Engineering's
$\Pi_1$  significance for
traditional arithmetics.)

\end{document}